\documentclass{article}\begin{document}
\centerline{\bf\Large Generalized Integrals and Solvability}
\vskip .4in   

\centerline{N.D. Bagis}
\centerline{Department of Informatics, Aristotele University}
\centerline{Thessaloniki, Greece}
\centerline{nikosbagis@hotmail.gr}
\vskip .2in

\[
\]

\textbf{Keywords}: Generalized Integrals; Lagrange Equation; Evaluations; Solvability.
\[
\]
\centerline{\bf Abstract}

\begin{quote}

Based on previous work we consturct an equation (Lagrange equation) and relate it with a system of  generalized integrals and differential equations in such a way to provide useful evaluations and connections between them. 

\end{quote}

\section{The inversion problem in the complex analog}

The Lagrange inversion formula states: If $f(A)$ is analytic in a disck $D\subset \textbf{C}$ of center zero and $f(A)\neq 0$ in $D$, then in some neighborhood arroud $0$ the equation
\begin{equation}
\frac{w}{f(w)}=q
\end{equation}
have solution 
\begin{equation}
w=w(q)=\sum^{\infty}_{n=1}c_nq^n.
\end{equation}
The coefficients $c_n$ are given from
\begin{equation}
c_n=\frac{1}{n!}\left[\left(\frac{D}{Dh}\right)^{n-1}(f(h)^n)\right]_{h=0}\textrm{, }n=1,2,\ldots.
\end{equation}
Moreover we can extend the above theorem (formula) to
\begin{equation}
g_0(w(q))=g_0(0)+\sum^{\infty}_{n=1}\frac{q^n}{n!}\left[\left(\frac{D}{Dh}\right)^{n-1}\left(g_0'(h)(f(h))^{n}\right)\right]_{h=0},
\end{equation}
where $g_0(A)$ is also analytic. Setting $g_0(A)=e^A$, we get
\begin{equation}
e^{w(q)}=1+\sum^{\infty}_{n=1}\frac{q^n}{n!}\left[\left(\frac{D}{Dh}\right)^{n-1}\left(e^h(f(h))^{n}\right)\right]_{h=0}
\end{equation}
Also in view of [3] Theorem 20, it holds the following formula
\begin{equation}
e^{w(q)}=\prod^{\infty}_{n=1}\left(1-q^n\right)^{-\frac{1}{n}\sum_{d|n}\frac{\mu\left(n/d\right)}{\Gamma(d)}\left[\left(\frac{D}{Dh}\right)^{d-1}(f(h)^d)\right]_{h=0}}.
\end{equation} 
Here $\mu(n)$ is the Moebius $\mu$ function.\\

Set now
\begin{equation}
a_n=nc_n.
\end{equation}
Then in view of [2], in the function
\begin{equation}
\frac{1}{P(z)}=\sum^{\infty}_{n=1}a_nq^n\textrm{, }q=e(z):=e^{2\pi i z}\textrm{, }Im(z)>0,
\end{equation}
is attached a differential equation
\begin{equation}
X'(A)+2^{4/3}A^{-2/3}\left(1-A^2\right)^{-1/3}P\left(X(A)\right)=0.
\end{equation}
If $m^{*}(z)$ is the elliptic singular modulus defined as (see [4]):
\begin{equation}
m^{*}(z):=\left(\frac{\theta_2\left(e^{i\pi z}\right)}{\theta_3\left(e^{i\pi z}\right)}\right)^2\textrm{, }Im(z)>0,
\end{equation}
where ($|q|<1$):
\begin{equation}
\theta_2(q)=\sum^{\infty}_{n=-\infty}q^{(n+1/2)^2}\textrm{, }\theta_3(q)=\sum^{\infty}_{n=-\infty}q^{n^2},
\end{equation}
are the ''null'' Jacobi theta functions. Then also in view of [3] we have that the function
\begin{equation}
Y(z)=X\left(m^{*}(2z)\right)
\end{equation}
satisfies 
\begin{equation}
Y'(z)=-4\pi i\cdot \eta(z)^4P\left(Y(z)\right),
\end{equation}
where $\eta(z)$ is the Dedekind eta function i.e.
$$
\eta(z)=q^{1/24}\prod^{\infty}_{n=1}(1-q^n)\textrm{, }q=e^{2\pi i z}\textrm{, }Im(z)>0.
$$
Moreover if 
\begin{equation}
F(z):=\int^{z}_{Y(i\infty)}\frac{dt}{P(t)},
\end{equation}
then
\begin{equation}
F\left(Y(z)\right)=-\sqrt[3]{2}B_0\left(m^{*}\left(2z\right)^2,\frac{1}{6},\frac{2}{3}\right),
\end{equation}
(here $B_0(z;a,b):=\int^{z}_{0}t^{a-1}(1-t)^{b-1}dt$ is the incomplete Beta function) and
\begin{equation}
F\left(Y\left(\frac{-1}{z}\right)\right)+F\left(Y(z)\right)=-\frac{\sqrt{3}\Gamma\left(\frac{1}{3}\right)^3}{\pi\sqrt[3]{2}}.
\end{equation}
\textbf{Note.} If we assume that, given a function $P(z)$ with $I=\int^{z}_{c}\frac{dt}{P(t)}$ known (i.e. for example elliptic integral or hyperelliptic integral), then as we present below, there exists a system of equations and relations between them. For example the function $Y(z)$ defined in (12) have modular properties. The integral $I$ itself is solution of an equation $w/f(w)=q$ and many other.\\ 
Also (see [3]):
$$
\exp\left(2\pi i\int^{z}_{Y(i\infty)}\frac{dt}{P(t)}\right)=\exp\left(2\pi i \int^{z}_{i\infty}\frac{dt}{P(t)}\right)\exp\left(2\pi i\int^{i\infty}_{Y(i\infty)}\frac{dt}{P(t)}\right)=
$$
$$
=\exp\left(w(q)\right)q^{C_0}=
$$
$$
=\prod^{\infty}_{n=1}\left(1-q^n\right)^{-\frac{1}{n}\sum_{d|n}a_d\mu(n/d)}
=\exp\left(8\pi^2\int^{Y^{(-1)}(z)}_{i\infty}\eta(t)^4dt\right).\eqno{(16.1)}
$$
Note that if $Im(z_1),Im(z_2)>0$, then
\begin{equation}
2\pi i\int^{z_2}_{z_1}\eta(t)^4dt=\left[\frac{1}{\sqrt[3]{4}}B_0\left(m^{*}(2t)^2;\frac{1}{6};\frac{2}{3}\right)\right]^{t=z_2}_{t=z_1},
\end{equation}
where $\eta(z)$ is the Dedekind's eta function
\begin{equation}
\eta(z)=q^{1/24}\prod^{\infty}_{n=1}\left(1-q^n\right)\textrm{, }q=e(z)\textrm{, }Im(z)>0.
\end{equation}
Hence:\\ 
\\
\textbf{Theorem 0.0}\\
Assume $q=e(z)$, $Im(z)>0$, then
$$
q^{C_0}\exp\left(w(q)\right)=
$$
$$
=\prod^{\infty}_{n=1}\left(1-q^n\right)^{-\frac{1}{n}\sum_{d|n}a_d\mu(n/d)}
=\exp\left(8\pi^2\int^{Y^{(-1)}(z)}_{i\infty}\eta(t)^4dt\right).\eqno{(18.1)}
$$
\\
\textbf{Theorem 0.1}\\
Assuming the general product
$$
\prod^{\infty}_{n=1}\left(1-e\left(nz\right)\right)^{-\chi(n)}\textrm{, }e(z):=e^{2\pi i z}\textrm{, }Im(z)>0,
$$
where 
$$
\chi(n)=\frac{1}{n}\sum_{d|n}a_d\mu\left(\frac{n}{d}\right),\eqno{(18.2)}
$$
we get the modular relation
$$
\prod^{\infty}_{n=1}\left(1-e\left(nY(z)\right)\right)^{\chi(n)}\prod^{\infty}_{n=1}\left(1-e\left(nY\left(\frac{-1}{z}\right)\right)\right)^{\chi(n)}=e^{i 4^{1/3}\sqrt{3}\Gamma\left(\frac{1}{3}\right)^3},\eqno{(18.3)}
$$
where $Y(z)$ is that of (12),(13).\\
\\

Now according to conection (7) we have
\begin{equation}
2\pi i\int^{z}_{i\infty}\frac{dt}{P(t)}=w(q).
\end{equation}
Hence
\begin{equation}
P(z)=\frac{1}{qw'(q)}.
\end{equation}
From (15) then we get
\begin{equation}
w\left(e^{2\pi i Y(z)}\right)=-2\pi i\sqrt[3]{2}B_0\left(m^{*}(2z)^2;\frac{1}{6};\frac{2}{3}\right)+c.
\end{equation}
Hence
\begin{equation}
w\left(e^{2\pi i X(A)}\right)=-2\pi i \sqrt[3]{2}B_0\left(A^2;\frac{1}{6};\frac{2}{3}\right)+c.
\end{equation}
But it is known that (see [3] Theorem 12)
\begin{equation}
X(A)=h\left(\sqrt[3]{2}B_0\left(A^2;\frac{1}{6};\frac{2}{3}\right)\right)
\end{equation}
Hence easily
\begin{equation}
e^{2\pi i h(A)}=w^{(-1)}\left(-2\pi i A+c\right)
\end{equation}
and 
\begin{equation}
e^{2\pi i h(A)}=\frac{-2\pi i A+c}{f\left(-2\pi i A+c\right)}.
\end{equation}
But relation (24) gives (we use the notation $q_A=e(A)$):
$$
w\left(e^{2\pi i A}\right)=-2\pi i h_i(A)+c\Rightarrow
$$
\begin{equation}
h_i(A)=-\frac{w\left(q_A\right)}{2\pi i}+\frac{c}{2\pi i}
\end{equation}
and
\begin{equation}
h_i'(A)=-w'(q_A)q_A\Leftrightarrow h_i'(A)P(A)=-1
\end{equation}
Hence we get also that 
$$
h_i\left(i\infty\right)=\frac{c}{2\pi i}\Rightarrow h\left(\frac{c}{2\pi i}\right)=i\infty\eqno{(28.1)}
$$
and
$$
Y(i\infty)=X(0).\eqno{(28.2)}
$$
From the analysis given in [3] we have (Corollary 1, eq 33) we have
\begin{equation}
5\int^{y(A)}_{0}\frac{dt}{t\sqrt[6]{t^{-5}-11-t^5}}=-\frac{w\left(e^{2\pi i A}\right)}{2\pi i}+\frac{c}{2\pi i}.
\end{equation}
However if we introduce the function $F_1(A)$ (as in [1]) such that
\begin{equation}
F_1'(A)=5^{-1}F_1(A)\sqrt[6]{F_1(A)^{-5}-11-F_1(A)^5},
\end{equation}
then
\begin{equation}
F_1^{(-1)}(A)=6A^{5/6}F_{Ap}\left[\frac{1}{6},\frac{1}{6},\frac{1}{6};\frac{7}{6};\frac{-2A^5}{11+5\sqrt{5}},\frac{-2A^5}{11-5\sqrt{5}}\right],
\end{equation}
where 
\begin{equation}
F_{Ap}(a,b_1,b_2;c;x,y)=\sum^{\infty}_{m,n=0}\frac{(a)_{m+n}(b_1)_m(b_2)_n}{(c)_{m+n}m!n!}x^my^n\textrm{, }|x|<1\textrm{, }|y|<1,
\end{equation}
is the first Appell function and
\begin{equation}
5\int^{F_1(A)}_{0}\frac{dt}{t\sqrt[6]{t^{-5}-11-t^5}}=A.
\end{equation}
\\
\textbf{Theorem 1.01} (see [6],[3]).\\
We assume that for a given $\chi(n)$, the arithmetical function $a_n$ is such that
\begin{equation}
q^{-C_0}\prod^{\infty}_{n=1}\left(1-q^n\right)^{-\chi(n)}=\textrm{Algebraic Number, }\forall q=e(z)\textrm{, }z\in D,
\end{equation}
where 
$$
D=\{z\in \textbf{C}: z=r_1+i\sqrt{r_2}, r_1, r_2\in \textbf{Q}\textrm{ and }r_2>0\}.
$$ 
Then
$$
q^{-C_0}\prod^{\infty}_{n=1}\left(1-q^n\right)^{-\chi(n)}=S_1(z)\textrm{, }q=e(z)\textrm{, }Im(z)>0
$$
and
$$
S_1(z)=S\left(m^{*}(2z)\right),
$$
where $S(z)$ is also algebraic function (from the set of algebraic numbers to the set of algebraic numbers) and 
$$
S_1(z)=c_1e^{w(q)}\textrm{, }q=e(z)\textrm{, }Im(z)>0.
$$
The function $w(q)$ can evaluated from the $a_n$ and relations (8),(7),(2). Moreover if for all $z\in D$, the function $S_1(z)$ is algebraic and $S(A)$ is also algebraic, then these two functions can constructed explicity using lower degree algebraic numbers and (hence) find their complete form (see [3],[6]).\\
\\
\textbf{Definition.}\\
A rational valued arithmetic function $X(n)$ will be called ''Moebius Periodic'' if exist $T$ odd positive integer (for even periods $T$ hold similar arguments) such that:\\
$\textbf{I.}$ $X(n)$ is $T-$periodic i.e. for all $n,k\in\{0,1,2,\ldots\}$, we have $\chi(k):=X(k+n T)=X(k)$\\and\\
$\textbf{II.}$ In every period interval hold $a_0=a_T=0$, $a_1=a_{T-1}$, $a_2=a_{T-2}$, $\ldots$, $a_{(T-1)/2}=a_{(T+1)/2}$.\\
\\
\textbf{Remarks.}\\
It is easy to see someone that if $X(n)$ is a periodic arithmetical function, with period $T$ positive integer, then
$$
\prod^{\infty}_{n=1}\left(1-q^n\right)^{X(n)}=\prod^{T}_{k=1}\left(q^k;q^T\right)^{X(k)}_{\infty},
$$ 
where 
$$
\left(a;q\right)_{\infty}:=\prod^{\infty}_{n=0}\left(1-aq^n\right).
$$
If also $X(n)$ is even, then $X(T-n)=X(-n)=X(n)$. Hence
$$
\left(\prod^{\infty}_{n=1}\left(1-q^n\right)^{X(n)}\right)^2=\left
(\prod^{\infty}_{n=1}\left(1-q^n\right)^{X(n)}\right)\left(\prod^{\infty}_{n=1}\left(1-q^n\right)^{X(n)}\right)=
$$
$$
=\left(\prod^{T}_{k=1}\left(q^k;q^T\right)^{X(k)}_{\infty}\right)\left(\prod^{T}_{k=1}\left(q^k;q^T\right)^{X(k)}_{\infty}\right)=
$$
$$
=\left(\prod^{T}_{k=1}\left(q^k;q^T\right)^{X(k)}_{\infty}\right)\left(\prod^{T}_{k=1}\left(q^{k-T};q^T\right)^{X(T-k)}_{\infty}\right)=
$$
$$
=\left(\prod^{T}_{k=1}\left(q^k;q^T\right)^{X(k)}_{\infty}\right)\left(\prod^{T}_{k=1}\left(q^{T-k};q^T\right)^{X(k)}_{\infty}\right)=
$$
$$
=\prod^{T}_{k=1}\left(\left(q^k;q^T\right)_{\infty}\left(q^{T-k};q^T\right)_{\infty}\right)^{X(k)}=
$$
$$
=q^{T/24\sum^{T}_{k=1}X(k)}\eta(Tz)^{-\sum^{T}_{k=1}X(k)}\prod^{T}_{k=1}\theta\left(\frac{T}{2},\frac{T}{2}-k;q\right)^{X(k)}.
$$
Hence we get the next\\
\\
\textbf{Theorem 1.02}\\
If $X(n)$ is $T-$periodic and even arithmetical function, then
$$
\left(\prod^{\infty}_{n=1}\left(1-q^n\right)^{X(n)}\right)^2=q^{T/24\sum^{T}_{k=1}X(k)}\eta(Tz)^{-\sum^{T}_{k=1}X(k)}\prod^{T}_{k=1}\theta\left(\frac{T}{2},\frac{T}{2}-k;q\right)^{X(k)},
$$
where $q=e(z)$, $Im(z)>0$ and 
$$
\theta(a,b;q):=\sum^{\infty}_{n=-\infty}(-1)^nq^{an^2+bn},
$$
$a,b$ rationals, $a>0$.\\
\\

On the same basis one has\\
\\
\textbf{Theorem 1.03} (see [6]).\\
Let $g(q)$ be a function analytic in $(-1,1)$ and set
$$
X(n)=\frac{1}{n}\sum_{d|n}\frac{g^{(d)}(0)}{\Gamma(d)}\mu\left(\frac{n}{d}\right).
$$
Define
$$
C=\sum^{\left[\frac{T-1}{2}\right]}_{j=1}\left(-\frac{j}{2}+\frac{j^2}{2T}+\frac{T}{12}\right)X(j).
$$
Then if $X(n)$ is $T-$Moebius periodic we have:\\
\textbf{I.} 
$$
q^{-C}e^{g(q)}=\textrm{Algebraic Number}
$$
at the points $q=e^{2\pi i z}$, where $z\in D$.\\
\textbf{II.}\\
The function $e^{g(q)}$ can be expanded in finite terms of $\theta-$theta functions and the Dedekind eta function $\eta(z)$:
$$
e^{g(q)}=\eta(Tz)^{\sum^{\left[\frac{T-1}{2}\right]}_{j=1}X(j)}\prod^{\left[\frac{T-1}{2}\right]}_{j=1}\theta\left(\frac{T}{2},\frac{T}{2}-j;q\right)^{-X(j)},
$$ 
where 
$$
\theta(a,b;q)=\sum^{\infty}_{n=-\infty}(-1)^n q^{an^2+b n}\textrm{, }q=e(z)\textrm{, }Im(z)>0
$$
and $a,b$ are rationals, with $a>0$.\\
\\ 
\textbf{Theorem 1.04}\\
Assume the product
$$
\Pi=q^{-C_0}\prod^{\infty}_{n=1}\left(1-q^n\right)^{-1/n\sum_{d|n}a_d\mu(n/d)}\textrm{, }|q|<1,
$$
where $a_n:\textbf{N}\rightarrow\textbf{C}$ is $T-$periodic arithmetical function with $T$ positive integer. Then
$$
\Pi=\exp\left(w(q)\right)=\exp\left(\sum^{\infty}_{n=1}a_n\frac{q^n}{n}\right)=\exp\left(\sum^{T}_{n=1}a_n\sum^{\infty}_{j=0}\frac{q^{n+j T}}{n+j T}\right)\textrm{, }|q|<1, 
$$
which is a finite function of the form
$$
\exp\left(w(q)\right)=\prod^{N_1}_{j=0}\left(1-\zeta(j,T;a_n)q\right)^{c(j,T;a_n)}.
$$
The $\zeta(j,T;a_n)$, $c(j,T;a_n)$, for $j=0,1,\ldots,N_1$ being algebraic functions, when $a_n$ are rationals.\\
\\
\textbf{Theorem 1.05}\\
Assume that $a_n$ is any arithmetical function. Then Taylor series
$$
\sum^{\infty}_{n=1}\frac{a_n}{n}q^n,
$$
is known, iff the Borcherds product
$$
q^{-C_0}\prod^{\infty}_{n=1}\left(1-q^n\right)^{-1/n\sum_{d|n}a_d\mu(n/d)},
$$
is known.\\
\textbf{Remark.}\\
By this way the evaluation of a Borcherds products are such easy as evaluating a Taylor series.\\
\\
\textbf{Theorem 1.1}\\
The functions $y(A)$ and $w(q_A)$ are related with the following identity:
\begin{equation}
F_1\left(-\frac{w\left(q_A\right)}{2\pi i}+\frac{c}{2\pi i}\right)=y(A).
\end{equation}
Also holds (see Theorem 22.1 relation (239) below):
$$
y(Y(A))=R(q_A),\eqno{(34.1)}
$$
where $q_A=e(A)$, $Im(A)>0$ and 
$$
R(q_A)=v(A)=q_{A}^{1/5}\prod^{\infty}_{n=1}\left(1-q_{A}^{n}\right)^{(n|5)},
$$ 
is the Rogers-Ramanujan continued fraction. Also
$$
y_i\left(A\right)=Y\left(v_i(A)\right).
$$
\\

Moreover we have
\begin{equation}
\frac{1}{G(y(A))}=h_i'(A)=-w'\left(q_A\right)q_A.
\end{equation}
Hence
\begin{equation}
-w'\left(q_A\right)q_AG\left(F_1\left(-\frac{w\left(q_A\right)}{2\pi i}+\frac{c}{2\pi i}\right)\right)=1\Leftrightarrow
\end{equation}
\begin{equation}
G\left(y(A)\right)=-P(A).
\end{equation}
Continuing our arguments we have
\begin{equation}
G\left(F_1\left(-\int^{P_i(A)}_{i\infty}\frac{dt}{P(t)}+\frac{c}{2\pi i}\right)\right)=-A.
\end{equation}
This is true because
$$
-\frac{w\left(q_A\right)}{2\pi i}=-\int^{A}_{i\infty}\frac{dt}{P(t)}.
$$
Hence (38) gives
$$
-\int^{P_i(A)}_{i\infty}\frac{dt}{P(t)}+\frac{c}{2\pi i}=F_1^{(-1)}\left(G^{(-1)}(-A)\right)\Rightarrow
$$
$$
-\frac{P_i'(A)}{A}=\frac{d}{dA}\left(F_1^{(-1)}\left(G^{(-1)}(-A)\right)\right)\Rightarrow
$$
\begin{equation}
{P^{(-1)}}{'}(A)=-A\frac{d}{dA}\left(F_1^{(-1)}\left(G^{(-1)}(-A)\right)\right).
\end{equation}
Also
$$
-{P^{(-1)}}{'}\left(e(A)\right)2\pi i
=-e(A)2\pi i{F_1^{(-1)}}{'}\left(G^{(-1)}\left(-e(A)\right)\right){G^{(-1)}}{'}\left(-e(A)\right)\Leftrightarrow
$$
$$
-2\pi i \int {P^{(-1)}}{'}\left(e(A)\right)dA=F_1^{(-1)}\left(G^{(-1)}\left(-e(A)\right)\right).
$$
Hence assuming that ${P^{(-1)}}{'}(A)=H'(A) A$, we get 
$$
H(A)=-F_1^{(-1)}\left(G^{(-1)}\left(-A\right)\right)+c\Leftrightarrow c-H(A)=F_1^{(-1)}\left(G^{(-1)}(-A)\right)\Leftrightarrow
$$
$$
G\left(F_1\left(c-H(A)\right)\right)=-A.
$$
However we have the next theorems.\\
\\
\textbf{Theorem 2.}\\
Assuming $q=e^{2\pi i A}$, $Im(A)>0$, we have
\begin{equation}
\int^{-w\left(q\right)/(2\pi i)+c/(2\pi i)}_{c'^{*}}G(F_1(t))dt=A.
\end{equation}
Hence given $G$ we can find $w(q)$ and the oposite.\\
\\
\textbf{Proof.}\\
Integrate (36).\\
\\
\textbf{Theorem 2.1}\\
It is easy to see someone that holds in general
$$
h'(A)=G\left(F_1(A)\right)\eqno{(40.1)}
$$
and
$$
y(A)=F_1\left(h_i(A)\right).\eqno{(40.2)}
$$
\\
\textbf{Theorem 3.}\\
Knowing $G(A)$ (resp. $w(q_A)$) we can find $w(q_A)$ (resp. $G(A)$) from Theorem 2 and then $P(A)$ from the identity
\begin{equation}
G\left(F_1\left(-\frac{w\left(q_A\right)}{2\pi i}+\frac{c}{2\pi i}\right)\right)=-P(A).
\end{equation}
Also there exists the relations 
\begin{equation}
w(q_A)=2\pi i\int^{A}_{i\infty}\frac{dt}{P(t)}\textrm{ and }P(A)=\frac{1}{q_Aw'(q_A)},
\end{equation}
where $q_A=e^{2\pi i A}$, $Im(A)>0$.\\
\\
\textbf{Theorem 4.}\\
Given the function $G(A)$ and assuming function $y(A)$ is solution to the problem
\begin{equation}
5\int^{y(A)}_{0}\frac{G(t)}{t\sqrt[6]{t^{-5}-11-t^5}}dt=A,
\end{equation}
then $P(A)$ is such that
\begin{equation}
{P^{(-1)}}{'}(A)=-A\frac{d}{dA}\left(F_1^{(-1)}\left(G^{(-1)}(-A)\right)\right)
\end{equation}
and $y(A)$ is solution of the semialgebraic equation
\begin{equation}
G(y(A))+P(A)=0.
\end{equation}
The function $F_1(A)$ is the known function defined in (30),(31),(33). Also $P(A)$ is given from
\begin{equation}
P(A)=\frac{1}{q_Aw'(q_A)}\textrm{, }q_A=e(A).
\end{equation}
\\
\textbf{Theorem 5.}\\
If $P(A)$ is the function (8), then the function
\begin{equation}
w(q)=2\pi i\int^{z}_{i\infty}\frac{dt}{P(t)}\textrm{, }q=e(z)\textrm{, }Im(z)>0,
\end{equation} 
is the solution of
\begin{equation}
\frac{w(q)}{f(w(q))}=q,
\end{equation}
where $f(A)$ is given from
\begin{equation}
\frac{f'(A)}{f(A)}=\frac{1}{A}+G\left(F_1\left(-\frac{A}{2\pi i}+\frac{c}{2\pi i}\right)\right).
\end{equation}
Also then
\begin{equation}
y(A)=F_1\left(-\frac{w\left(q_A\right)}{2\pi i}+\frac{c}{2\pi i}\right).
\end{equation}
\\
\textbf{Theorem 6.}\\
Assume that exists function $P_0(A)$ such that
\begin{equation}
G\left(F_1(A)\right)=-\frac{1}{c-2\pi i A}+P_0\left(c-2\pi i A\right).
\end{equation}
Then
\begin{equation}
y(A)=F_1\left(-\frac{w\left(e^{2\pi i A}\right)}{2\pi i}+\frac{c}{2\pi i}\right),
\end{equation}
where $w(q)$ is solution of the equation
\begin{equation}
w(q)\exp\left(C-\int^{w(q)}_{c'}P_0(t)dt\right)=q.
\end{equation}
\textbf{Remarks.}\\ 
\textbf{1)} Given a $G(A)$ of the form (51) we evaluate $w(q)$ from (53). Also we find $x(A)=w(e(P_i(A)))$ from  
\begin{equation}
A+P_{0}(x(A))=\frac{1}{x(A)}
\end{equation}
and $P(A)$ from 
\begin{equation}
P(A)=\frac{1}{w(q)}-P_0(w(q)).
\end{equation}
Such $G$ as in (51) is called ''proper''.\\
\textbf{2)} Also holds
\begin{equation}
A+G\left(F_1\left(\frac{c-x(A)}{2\pi i}\right)\right)=0
\end{equation}
and
\begin{equation}
x'(A)=2\pi i\frac{P_i'(A)}{A}\textrm{, }x_i(A)=\frac{w_i'(A)}{w_i(A)}.
\end{equation}
The function $V(A)$ defined as
\begin{equation}
V(x_i(A))=\frac{A^4P_0'(A)}{1+A^2P_0'(A)},
\end{equation}
satisfies the equations
\begin{equation}
V(A)=\frac{d}{dA}\left(\int^{x(A)}_{c}P_0'(t)t^2dt\right)
\end{equation}
and
\begin{equation}
x'(A)+x(A)^2=V(A).
\end{equation}
\\

Set now (see [1],[2]):
$$
C(\nu):=\int^{1}_{0}f(-q)^4q^{-5/6}R(q)^{5\nu}dq=
$$
$$
=\Gamma\left(\frac{5}{6}\right)\left(\frac{11+5 \sqrt{5}}{2}\right)^{-\frac{1}{6}-\nu} \frac{\Gamma\left(\frac{1}{6}+\nu\right)}{\Gamma(1+\nu)} {}_2F_1\left(\frac{1}{6},\frac{1}{6}+\nu;1+\nu;\frac{11-5 \sqrt{5}}{11+5\sqrt{5}}\right),
$$
where $\nu\geq0$.\\
\\
\textbf{Theorem 6.1} (see [2] Theorem 18.1)\\
Assume the function $G(t)$ is given near the origin by (under certain converging conditions):
\begin{equation}
G(t)=\sum^{\infty}_{m=0}a_mt^{p_m},
\end{equation}
where $p_m$ is any increasing sequence of positive real numbers with $\lim p_m=+\infty$. Then the function $y(A)$ defined as 
\begin{equation}
5\int^{y(A)}_{0}\frac{G(t)}{t\sqrt[6]{t^{-5}-11-t^5}}dt=A,
\end{equation}
have the following property: Every equation of the form
\begin{equation}
y(x)=\alpha \frac{\sqrt{5}-1}{2},
\end{equation}
have solution $x$ such that
\begin{equation}
x=\sum^{\infty}_{m=0}a_m^{*}(\alpha)C\left(5^{-1}p_m\right),
\end{equation}
where $a^{*}_m(\alpha)$ is such that
\begin{equation}
G(x\alpha)\sqrt[6]{\frac{x^{-5}-11-x^5}{(\alpha x)^{-5}-11-(\alpha x)^5}}=\sum^{\infty}_{m=0}a^{*}_{m}(\alpha)x^{p_m}.
\end{equation}
Hence if we set $\xi^{-1}=\frac{\sqrt{5}-1}{2}$, then the inverse of $y(A)$ is
\begin{equation}
y_i(A)=\sum^{\infty}_{m=0}a^{*}_m(\xi A)C\left(5^{-1}p_m\right),
\end{equation}
provited the convergence of (61),(62),(64),(65),(66).\\
\\
\textbf{Theorem 6.1.1}\\
It holds (see [2]: Note equation (69))
$$
h\left(\frac{1}{\sqrt[3]{4}}B_0\left(A^2;\frac{1}{6},\frac{2}{3}\right)\right)=s_i(A)\eqno{(66.1)}
$$
and
$$
h_i(A)=-\frac{w(q_A)}{2\pi i}+\frac{c}{2\pi i}.\eqno{(66.2)}
$$
\textbf{Remarks.}\\
\textbf{i)} We often use the notation $f_i(A)$ to stand for $f^{(-1)}(A)$, i.e. the inverse function of $f(A)$.\\
\textbf{ii)} When equation (1) is given, then we can find from (66.2) the function $h_i(A)$ and hence from (66.1) the function $s_i(A)$. However from Theorem 10 of [2], when the function $G(t)$ is written as
$$
G(t)=\frac{G_0(t)}{\sigma(F_i(t))},\eqno{(66.3)}
$$
where $F(A)=R\left(e^{\pi i m^{*}_i(A)}\right)$, $m^{*}(A)$ is the eliptic singular modulus (see Theorem 22.1 below) and 
$$
m_i^{*}(A)=i\frac{K\left(\sqrt{1-A^2}\right)}{K(A)}\textrm{, }K(A)=\frac{\pi}{2}\cdot {}_2F_1\left(\frac{1}{2},\frac{1}{2};1;A^2\right).\eqno{(66.4)}
$$
Then 
$$
y(A)=R\left(e^{\pi i m_i^{*}(s(A))}\right)\textrm{ and }s'(A)=\sigma(s(A)).\eqno{(66.5)}
$$
The function $G_0(A)$ is
$$
\frac{3\cdot 2^{2/3}\sqrt{G_0(A)}}{\sqrt[6]{1+\sqrt{1-8 G_0(A)^3}}}\cdot {}_2F_1\left[\frac{1}{6},\frac{1}{3};\frac{7}{6};\frac{1}{2}\left(1-\sqrt{1-8 G_0(A)^3}\right)\right]=
$$
$$
=6A^{5/6} F_{Ap}\left[\frac{1}{6};\frac{1}{6},\frac{1}{6};\frac{7}{6};\frac{11-5 \sqrt{5}}{2}A^5,\frac{11+5 \sqrt{5}}{2}A^5\right]=F^{(-1)}_1(A).\eqno{(66.6)}
$$
\textbf{iii)} More general if
$$
G(x)=\sum^{\infty}_{n=0}G_nx^{p_n},\eqno{(66.6.1)}
$$
where $p_n$ is increasing sequence of positive real numbers, with $\lim p_n=+\infty$, or $p_n$ strictly increasing in $n\in\{0,1,\ldots,N\}$, with $p_{n}=0,\forall n>N$, then (see [2] page 27)
$$
h(A)=5F_1(A)^{5/6}\times
$$
$$
\times\sum^{\infty}_{n=0}G_n\frac{F_1(A)^{p_n}}{5/6+p_n}F_{Ap}\left[\frac{1}{6}+\frac{p_n}{5},\frac{1}{6},\frac{1}{6},\frac{7}{6}+\frac{p_n}{5},\rho_1 F_1(A)^5,\rho_2 F_1(A)^5\right],\eqno{(66.6.2)}
$$
where $\rho_1=\frac{11-5\sqrt{5}}{2}$, $\rho_2=\frac{11+5\sqrt{5}}{2}$. Also then
$$
y\left(5A^{5/6}\sum^{\infty}_{n=0}G_n\frac{A^{p_n}}{5/6+p_n}F_{Ap}\left[\frac{1}{6}+\frac{p_n}{5},\frac{1}{6},\frac{1}{6},\frac{7}{6}+\frac{p_n}{5},\rho_1 A^5,\rho_2 A^5\right]\right)=A.\eqno{(66.6.3)}
$$
\\
\textbf{Theorem 6.2}\\ 
\textbf{1)} Let $G$ be ''proper'' i.e.
\begin{equation}
G\left(F_1\left(\frac{c-A}{2\pi i}\right)\right)=-\frac{1}{A}+P_0(A),
\end{equation}
where $P_0(A)$ is some analytic function in a region. Then we can find $f(A)$ from the equation $P_0(A)=f'(A)/f(A)$. Solving for this $f$ the equation $\frac{w(q)}{f(w(q))}=q$, we find $w(q)$ as also from (52) function $y(A)$.\\However if we can not find $f(A)$ solving $f'(A)/f(A)=P_0(A)$, we can solve the equation (62) of Theorem 6.1 and find $y(A)$ using (66). Hence from relation (45) we can find $P(A)$. Knowing $P(A)$ and $P_0(A)$ we find $w(q)$ from equation (55) and hence $f(A)$. This can be used in Theorem 13.1 below.\\
\textbf{2)} Given a function $G$, there exists function $x(A)$ such that
\begin{equation}
A+G\left(F_1\left(\frac{c-x(A)}{2\pi i}\right)\right)=0.
\end{equation}
Then also exists $P_0(A)$ analytic such that
\begin{equation}
x^{(-1)}(A)=\frac{1}{A}-P_0(A).
\end{equation}
Given a potential $V(A)$, we solve 
$$
x'(A)+x(A)^2=V(A).
$$
Then $P_0(A)$ and $x(A)$ are related as
$$
x^{(-1)}(A)=\frac{1}{A}-P_0(A).
$$
For this $x(A)$ we have if $y_1(A)=\exp\left(\int^{A}_{c}x(t)dt\right)$, then $y_1(A)$ is solution of $y_1''(A)=V(A)y_1(A)$ and $G(A)$ is given from (67). Also we can find $w(q)$, $P(A)$ from Theorem 6 and relation (55) resp. Hence from  $G(y(A))=-P(A)$ again we can find $y(A)$. This is another way to find the solution $y(A)$ of the problem (62).\\
\\
\textbf{Theorem 6.2.0}\\
Assume the Ricatti equation
$$
x'(A)+x(A)^2=V(x),\eqno{(69.01)}
$$
where the potential $V(x)$ is $T-$periodic trigonometric polynomial:
$$
V(A)=\sum^{N}_{n=0}c_ne^{-2\pi i n A/T}\textrm{, }A\in \left[0,T\right].\eqno{(69.02)}
$$
Then setting $x(A)=\frac{Y'_1(A)}{Y_1(A)}$, we get
$$
Y''_1(A)=V(A)Y_1(A).\eqno{(69.03)}
$$
Also if $f_0(A)=\widehat{Y}_1(A)$, then
$$
\sum^{N}_{n=0}c_nf_0\left(A+\frac{2\pi n}{T}\right)=-A^2f_0(A)\eqno{(69.04)}
$$
and
$$
x(A)=\frac{d}{dA}\log\left(\frac{1}{2\pi}\int^{\infty}_{-\infty}f_0(\gamma)e^{i\gamma A}d\gamma\right).\eqno{(69.05)}
$$
\textbf{Remarks.} When $N\rightarrow\infty$, we have
$$
V(A)=\sum^{\infty}_{n=0}c_ne^{-2\pi i n A/T}\textrm{, }A\in \left[0,T\right],\eqno{(69.05)}
$$
$$
Y''(A)=V(A)Y(A)
$$
and
$$
Y(A)=\frac{1}{2\pi}\int^{\infty}_{-\infty}f_0(\gamma)e^{i\gamma A}d\gamma,\eqno{(69.05)}
$$
where
$$
\sum^{\infty}_{n=0}c_nf_0\left(A+\frac{2\pi n}{T}\right)=-A^2f_0(A).\eqno{(69.06)}
$$
This last equation is equivalent to (set $A\rightarrow \frac{2\pi}{T}A$ and $f_0\left(A\frac{2\pi }{T}\right)=f_1\left(A\right)$, $c_n\frac{T^2}{4\pi^2}=c_n^{*}$):
$$
\sum^{\infty}_{n=0}c^{*}_nf_1\left(A+n\right)=-A^2f_1(A),\eqno{(69.07)}
$$
where $c^{*}_n$ are given.\\
\\
\textbf{Proof.}\\
See [5].\\
\\

\[
\]

Hence $G(A)$ and $V(A)$ are related with the equation
$$
-\frac{d}{dA}\left(-G\left(F_1\left(\frac{c-A}{2\pi i}\right)\right)\right)=\left[A^2-V\left(-G\left(F_1\left(\frac{c-A}{2\pi i}\right)\right)\right)\right]^{-1}.\eqno{(69.1)}
$$
Also given the potential $V(A)$ we can find $P_0(A)$ from the equation
$$
V\left(\frac{1}{A}-P_0(A)\right)=\frac{A^4P_0'(A)}{1+A^2P_0'(A)},\eqno{(69.2)}
$$
which is simply the Ricatti equation 
$$
\psi^{(-1)}{'}(A)+\left(\psi^{(-1)}(A)\right)^2=V(A)\textrm{, where }\psi(A)=\frac{1}{A}-P_0(A).\eqno{(69.3)}
$$
Hence knowing $P_0(A)$, we find $w(q)$, $G(A)$ and $P(A)$ (from (53),(67) and (55)). From (45) we can find $y(A)$. In the oposite if we know $G(A)$, we can find $P_0(A)$ from (67) and hence $V(A)$. Also the equation (69.1) is a Ricatti equation.\\
\\
\textbf{Theorem 6.2.1}\\
Knowing $G(A)$, we can find $y(A)$ from (66). Hence from (45) we find $P(A)$. Using now (291) of Theorem 25.3 below, we find $f(A)$.\\ 
\\

Hence if we assume that $V(A)$ is involved, then exists $x_1(A)$ such that
\begin{equation}
A+V\left(F_1\left(\frac{c-x_1(A)}{2\pi i}\right)\right)=0\Rightarrow 
\end{equation}
\begin{equation}
x_1(A)=c-2\pi i F_1^{(-1)}\left(V^{(-1)}(-A)\right).
\end{equation}
Also from $x_1'(A)=2\pi i \frac{P_1^{(-1)}{'}(A)}{A}$ and $A+P_{10}(x_1(A))=\frac{1}{x_1(A)}$, we get
\begin{equation}
P_1^{(-1)}{'}(A)=AF_1^{(-1)}{'}\left(V_i(-A)\right)V_i'(-A)
\end{equation}
and
\begin{equation}
P_{10}(A)=\frac{1}{A}+V\left(F_1\left(\frac{c-A}{2\pi i}\right)\right).
\end{equation}
Also from
\begin{equation}
V_1(A)=2\pi i F_1^{(-1)}{'}\left(V_i\left(-A\right)\right)V_{i}'(-A)+\left(c-2\pi i F_1^{(-1)}\left(V_i(-A)\right)\right)^2,
\end{equation}
we get (Note that: $x_1'(A)+x_1(A)^2=V_1(A)$):
\begin{equation}
V_1\left(\frac{1}{A}-P_{10}(A)\right)=\frac{A^4P_{10}'(A)}{1+A^2 P_{10}'(A)}.
\end{equation}
Also
\begin{equation}
x(A)=c-2\pi i F_1^{(-1)}\left(G^{(-1)}(-A)\right)
\end{equation}
\begin{equation}
x_i(A)=-G\left(F_1\left(\frac{c-A}{2\pi i}\right)\right)
\end{equation}
and in the same way as (75) we get
\begin{equation}
V\left(\frac{1}{A}-P_0(A)\right)=\frac{A^4P_0'(A)}{1+A^2P_0'(A)}.
\end{equation}

We call distance of the functions $f,g$ the quantity
\begin{equation}
\textrm{dist}(f,g):=\left|f^{(-1)}(g(A))-g^{(-1)}(f(A))\right|.
\end{equation}
Then holds: $\textrm{dist}(f,f)=0$, $\textrm{dist}(f,g)=\textrm{dist}(g,f)$ and $\textrm{dist}(f,g)=0$ iff $h_1(x)=f^{(-1)}(g(A))$ is such that $h_1(h_1(A))=A$. Also if $\textrm{dist}(f,g)=0$ and $\textrm{dist}(g,h)=0$, then in general $\textrm{dist}(f,h)\neq 0$. Assume that $\textrm{dist}(f,g)=\textrm{dist}(g,h)=0$. Then exists $h_1(A),h_2(A)$, such that $f^{(-1)}(g(A))=g^{(-1)}(f(A))=h_1(A)$ and $g^{(-1)}(h(A))=h^{(-1)}(g(A))=h_2(A)$. Then $h_1(h_1(A))=A$ and $h_2(h_2(A))=A$. Also $f^{(-1)}(A)=h_1(g^{(-1)}(A))$ and $g^{(-1)}(A)=h_2(h^{(-1)}(A))$. Hence $f^{(-1)}(A)=h_1(h_2(h^{(-1)}(A)))$. Hence $f^{(-1)}(h(A))=h_1(h_2(A))\Rightarrow h^{(-1)}(f(A))=h_2(h_1(A))$. Hence 
$$
\textrm{dist}(f,h)=\left|h_1(h_2(A))-h_2(h_1(A))\right|.
$$

Define the conjugate of $\textrm{dist}$ with
$$
\textrm{dist}^{*}(f,g):=\left|f(g^{(-1)}(x))-g(f^{(-1)}(x))\right|.
$$
Then assuming for the curve  $C:(f(x),g(x))$, the change of variables $x\rightarrow k(x)$, we get $T(C):(f(k(x)),g(k(x)))$ and
$$
\textrm{dist}^{*}\left(f(k(x)),g(k(x))\right)=
$$
$$
=\left|f(k(k^{(-1)}(g^{(-1)}(x)))-g(k(k^{(-1)}(f^{(-1)}(x))))\right|=
$$
$$
=\left|f(g^{(-1)}(x))-g(f^{(-1)}(x))\right|=
$$
$$
=\textrm{dist}^{*}(f,g).
$$
Hence $\textrm{dist}^{*}$ is an invariant of the curve $C$. Also assume $\textrm{dist}^{*}(f,g)=0$, then $f(g^{(-1)}(x))=g(f^{(-1)}(x))$. Hence $f(g^{(-1)}(x))=\left(f(g^{(-1)}(x))\right)^{(-1)}$. Hence exists function $h(x)$ such that $f(g^{(-1)}(x))=h(x)$ and $h(h(x))=x$. Hence $f(x)=h(g(x))$. But then
$$
\textrm{dist}(f,g)=\left|f^{(-1)}(g(x))-g^{(-1)}(f(x))\right|=
$$
$$
=\left|f^{(-1)}(h^{(-1)}(f(x)))-f^{(-1)}(h(f(x)))\right|=0,
$$
since $h^{(-1)}(x)=h(x)$. Also we have:
$$
\textrm{dist}(f,g)=0\Leftrightarrow \textrm{dist}^{*}(f,g)=0.
$$

\[
\]

In our case we have 
$$
x_1^{(-1)}(A)=-V\left(F_1\left(\frac{c-A}{2\pi i}\right)\right)
$$
and
$$
x^{(-1)}(A)=-G\left(F_1\left(\frac{c-A}{2\pi i}\right)\right)\Rightarrow x(A)=c-2\pi i F_1^{(-1)}\left(G^{(-1)}(-A)\right). 
$$
Hence
\begin{equation}
x_1^{(-1)}(x(A))=-V\left(G_i(-A)\right)
\end{equation}
and similarly
\begin{equation}
x^{(-1)}(x_1(A))=-G\left(V^{(-1)}(-A)\right).
\end{equation}
Hence the distance of $x_1(A)$ and $x(A)$ is
$$
\textrm{dist}(x_1,x)=\left|x_1^{(-1)}(x(A))-x^{(-1)}(x_1(A))\right|=
$$
$$
=\left|-V\left(G_i(-A)\right)+G\left(V_i(-A)\right)\right|.
$$
From this last evaluation we have
\begin{equation}
\textrm{dist}\left(x_1(A),x(A)\right)=0\Leftrightarrow V\left(G_i\left(V(G_i(A)\right)\right)=A.
\end{equation}
Hence we get the next\\
\\
\textbf{Theorem 6.3}\\
Assume that $V,G$ are involved, then
$$
\textrm{dist}(x_1,x)=\textrm{dist}\left(G^{(-1)}(A),V^{(-1)}(A)\right).\eqno{(82.1)}
$$
Hence
\begin{equation}
\textrm{dist}(x_1,x)=0\Leftrightarrow \textrm{dist}(V,G)=0.
\end{equation}
\\

Assume now that $V(A)=h(G(A))$, with $(h\circ h)(A)=A$. Then $\textrm{dist}(x_1,x)=0$. Hence exists $h_1(A)$ such that $(h_1\circ h_1)(A)=A$ and $x(A)=h_1(x_1(A))$. Hence 
$$
x(A)=h_1\left(c-2\pi i F_1^{(-1)}\left(V^{(-1)}(-A)\right)\right).
$$
However if $G,V$ are real, continuous, invertible in a certain region, then $V\circ G_i$ is one to one. Hence the function $V(G_i(A))$ is strictly monotone. If $V(G_i(A))$ is strictly increasing, then $V(G_i(A))=A$ and thus $V(A)=G(A)\Rightarrow x_1(A)=x(A)$. But if $V(G_i(A))$ is strictly decreasing we have that exists a wide class of functions $V(G_i(A))$ such that (82) holds. Some examples are $1/A$, $c-A$, $\frac{-dA+b}{cA+d}$, $g_1\left(\frac{-dg_1^{(-1)}(A)+b}{cg_1^{(-1)}(A)+d}\right),\ldots$etc. Hence if
\begin{equation}
V\left(G_i(A)\right)=g_1\left(\frac{-dg_1^{(-1)}(A)+b}{cg_1^{(-1)}(A)+d}\right),
\end{equation} 
where $g_1$ is any smooth function, then (82) holds and
\begin{equation}
G\left(V_i(A)\right)=g_1\left(\frac{-dg_1^{(-1)}(A)+b}{cg_1^{(-1)}(A)+d}\right)\Rightarrow
\end{equation}
\begin{equation}
G(A)=g_1\left(\frac{-dg_1^{(-1)}(V(A))+b}{cg^{(-1)}_1(V(A))+d}\right).
\end{equation}
Bellow we will find lots of formulas that satisfy the functional equation $f(f(A))=A$. For now we state the next\\
\\
\textbf{Theorem 6.4}\\
Asume that $x_1(A)$, $x_2(A)$ are solutions of the equations $x_1'(A)+x_1(A)^2=V_1(A)$ and $x_2'(A)+x_2(A)^2=V_2(A)$ resp. and $V(A)$, $V_1(A)$ are defined from 
$$
A+V\left(F_1\left(\frac{c-x_1(A)}{2\pi i}\right)\right)=0\textrm{  and  }A+V_1\left(F_1\left(\frac{c-x_2(A)}{2\pi i}\right)\right)=0.\eqno{(86.1)}
$$
Then
$$
\textrm{dist}\left(x_1(A),x_2(A)\right)=\textrm{dist}\left(V^{(-1)}(A),V_1^{(-1)}(A)\right).\eqno{(86.2)}
$$
\\
\textbf{Theorem 7.}\\
\textbf{i)} It holds
\begin{equation}
F(z)=\frac{1}{2\pi i}w(q)+c_1\textrm{, }q=e(z)\textrm{, }Im(z)>0,
\end{equation}
where 
$$
c_1=\lim_{\sigma\rightarrow+\infty}\int^{+i\infty}_{Y(i\sigma)}\frac{dt}{P(t)}=\int^{+i\infty}_{Y(+i\infty)}\frac{dt}{P(t)}.\eqno{(87.1)}
$$
\textbf{ii)} If we define the function $g(A)$ to be 
\begin{equation}
g(A)=A\exp\left(-\int^{A}_{c}P_0(t)dt\right)
\end{equation}
and $P_0$ defined as in Theorem 6. Then the relation 
\begin{equation}
g^{(k)}(2\pi i (c_0-2c_1))=(-1)^k g^{(k)}(0)\textrm{, }c_0=-\frac{\sqrt{3}\Gamma\left(\frac{1}{3}\right)^3}{\pi\sqrt[3]{2}},
\end{equation}
is imposible when $g$ is not constant.\\
\textbf{Remarks.}\\
\textbf{i)} Condition (88) is equivalent to say that $g$ is analytic in $D_0$ and 
\begin{equation}
g\left(2 \pi i (c_0-2c_1)-z\right)=g(z)\textrm{, }\forall z\in D_0.
\end{equation}
\textbf{ii)} The set $D_0$ is subset of $\textbf{C}$ containg at least one circle with origin $0$ and radius greater than $2\pi |c_0-2c_1|>0$.\\
\\
\textbf{Proof.}\\
Assuming that (89) holds for every $k=0,1,2,\ldots$. We consider the Taylor series of $g$ arround $0$ and $2\pi i (c_0-2c_1)\neq 0$. We then have 
$$
e(Y(z))=w(e(Y(z)))\exp\left(-\int^{w(e(Y(z)))}_{c}P_0(t)dt\right)
$$
and
$$
e(Y(z))=\sum^{\infty}_{k=0}\frac{g^{(k)}(2\pi i (c_0-2c_1))}{k!}\left(w(e(Y(z)))-2\pi i (c_0-2c_1)\right)^k=
$$
$$
=\sum^{\infty}_{k=0}\frac{g^{(k)}(2\pi i(c_0-2c_1))}{k!}(2\pi i)^k\left(F(Y(z))-c_0+c_1\right)^k=
$$
$$
=\sum^{\infty}_{k=0}\frac{g^{(k)}(2\pi i(c_0-2c_1))}{k!}(2\pi i)^k\left(-F\left(Y\left(\frac{-1}{z}\right)\right)+c_1\right)^k=
$$
$$
=\sum^{\infty}_{k=0}\frac{g^{(k)}(2\pi i (c_0-2c_1))}{k!}(-1)^k\left(w\left(e\left(Y\left(\frac{-1}{z}\right)\right)\right)\right)^k=
$$
$$
=\sum^{\infty}_{k=0}\frac{g^{(k)}(0)}{k!}\left(w\left(e\left(Y\left(\frac{-1}{z}\right)\right)\right)\right)^k=e\left(Y\left(\frac{-1}{z}\right)\right).
$$
Hence $2\pi i Y\left(\frac{-1}{z}\right)=2\pi i Y\left(z\right)+2\pi i k_0$, $k_0\in\textbf{Z}$. Hence $Y\left(\frac{-1}{z}\right)=Y\left(z\right)+k_0$ and from realtion (16)
$$
F\left(z\right)+F\left(k_0+z\right)=c_0,
$$
for all $z\in D_0$. But from periodicity of $\frac{1}{P(z)}$ we have the existance of another constant $c_2$ such that $F(z+1)-F(z)=c_2$. Hence $F(z+k_0)-F(z)=k_0c_2$ and $F(z)+k_0 c_2+F(z)=c_0$. Hence $F$ is constant, which is imposible.\\
\\
\textbf{Example.}\\
Assume 
$$
G(A)=-\frac{1}{c-2\pi i F_1^{(-1)}(A)}+\cos\left(c-2\pi i F_1^{(-1)}(A)\right).
$$
Then $P_0(A)=\cos A$, and the function $y(A)$ such that
$$
5\int^{y(A)}_{0}\frac{G(t)}{t\sqrt[6]{t^{-5}-11-t^5}}=A
$$
is
$$
y(A)=F_1\left(\frac{c}{2\pi i}-\frac{w\left(e(A)\right)}{2\pi i}\right),
$$
where $w(q)$ is solution of
$$
w(q)e^{C_1-\sin(w(q))}=q.
$$
\\
\textbf{Theorem 8.}\\
Assume the function $h_0$ defined from the relations
\begin{equation}
g^{(k)}(2\pi i (c_0-2c_1))=(-1)^k g_2^{(k)}(0)\textrm{, }c_0=-\frac{\sqrt{3}\Gamma\left(\frac{1}{3}\right)^3}{\pi\sqrt[3]{2}},
\end{equation}
where
\begin{equation}
g_2(z)=h_0(g(z)).
\end{equation}
Then if $c_1$ denotes the constant $c_1=F(z)-\frac{w\left(e(z)\right)}{2\pi i}$ we have
\begin{equation}
g(2\pi i(c_0-2c_1)-z)=h_0(g(z)),
\end{equation}
\begin{equation}
w(A)+w(h_0(A))=c_2=2\pi i(c_0-2c_1)=\textrm{constant},
\end{equation}
\begin{equation}
e\left(Y\left(\frac{-1}{z}\right)\right)=h_0\left(e\left(Y\left(z\right)\right)\right)
\end{equation}
and
\begin{equation}
h_0(h_0(z))=z.
\end{equation}
$$
\textrm{dist}\left(Y\left(\frac{-1}{z}\right),Y\left(z\right)\right)=0.\eqno{(96.1)}
$$
\\
\textbf{Proof.}\\
Relation (95) can be shown as in Theorem 7. For to show (96) we have
$$
F(z)=\frac{w\left(e(z)\right)}{2\pi i}+c_1
$$
and
$$
F(Y(z))+F\left(Y\left(\frac{-1}{z}\right)\right)=c_0.
$$
Hence
$$
w\left(e\left(Y(z)\right)\right)+2\pi i c_1+w\left(e\left(Y\left(\frac{-1}{z}\right)\right)\right)+2\pi i c_1=2\pi i c_0\Leftrightarrow
$$
$$
w\left(e\left(Y(z)\right)\right)+w\left(h_0\left(e\left(Y(z)\right)\right)\right)=2\pi i (c_0-2c_1).
$$
Hence
$$
w\left(A\right)+w\left(h_0(A)\right)=2\pi i(c_0-2c_1).
$$
Setting where $h_0\rightarrow h_0^{(-1)}$ in the last equation, we have $w(h_0^{(-1)}(A))+w(A)=2\pi i (c_0-2c_1)$. Hence
$$
w\left(h_0^{(-1)}(A)\right)=w\left(h_0(A)\right)\Rightarrow h_0(h_0(A))=A.
$$
QED.\\
\\
\textbf{Theorem 9.}\\
We define the $B$ function to be such
\begin{equation}
h_0(A)=e(B(A)),
\end{equation}
and the $\lambda$ function  
\begin{equation}
\lambda(A)=B(e(A)).
\end{equation}
Then
\begin{equation}
\int^{A}_{i\infty}\frac{dt}{P(t)}+\int^{\lambda(A)}_{i\infty}\frac{dt}{P(t)}=c_0-2c_1,
\end{equation}
\begin{equation}
\lambda'(A)=-\frac{P(\lambda(A))}{P(A)},
\end{equation}
\begin{equation}
F(A)+F(\lambda(A))=\textrm{constant},
\end{equation}
\begin{equation}
F\left(\lambda(\lambda(A))\right)=F(A),
\end{equation}
where $F$,$\lambda$ are $1-$periodic
\begin{equation}
h_0(e(A))=e(\lambda(A)).
\end{equation}
There exists always integer $k=k(z)$ such that
\begin{equation}
Y\left(-\frac{1}{z}\right)=\lambda\left(Y(z)\right)+k
\end{equation}
and
\begin{equation}
2\pi i F(B(A))+w(A)=2\pi i (c_0-c_1).
\end{equation}
\\
\textbf{Proof.}\\
It holds
$$
w(e(A))+w\left(h_0(e(A))\right)=c_2\Leftrightarrow 
$$
$$
2\pi i \int^{A}_{i\infty}\frac{dt}{P(t)}+2\pi i \int^{B(e(A))}_{i\infty}\frac{dt}{P(t)}=c_2\Leftrightarrow
$$
$$
\int^{A}_{i\infty}\frac{dt}{P(t)}+\int^{\lambda(A)}_{i\infty}\frac{dt}{P(t)}=\frac{c_2}{2\pi i}.
$$
Setting $A\rightarrow \lambda(A)$, we get
$$
\int^{\lambda(\lambda(A))}_{i\infty}\frac{dt}{P(t)}=\frac{c_2}{2\pi i}-\int^{\lambda(A)}_{i\infty}\frac{dt}{P(t)}\Leftrightarrow
$$
$$
\int^{\lambda(\lambda(A))}_{i\infty}\frac{dt}{P(t)}=\int^{A}_{i\infty}\frac{dt}{P(t)}.
$$
Hence
$$
F\left(\lambda\left(\lambda(A)\right)\right)=F(A).
$$
For (101) we have
$$ 
Y\left(\frac{-1}{A}\right)=F^{(-1)}\left(c_0-F\left(Y(A)\right)\right)\Leftrightarrow
$$
$$
e\left(Y\left(\frac{-1}{A}\right)\right)=e\left(F^{(-1)}\left(c_0-F\left(Y(A)\right)\right)\right)\Leftrightarrow
$$
$$
h_0\left(e\left(Y(z)\right)\right)=e\left(F^{(-1)}\left(c_0-F\left(Y(A)\right)\right)\right)\Leftrightarrow
$$
$$
h_0\left(e\left(A\right)\right)=e\left(F^{(-1)}\left(c_0-F\left(A\right)\right)\right)\Leftrightarrow
$$
$$
e\left(B(e(A))\right)=e\left(F^{(-1)}\left(c_0-F\left(A\right)\right)\right)\Leftrightarrow
$$
$$
B\left(e(A)\right)=F^{(-1)}\left(c_0-F\left(A\right)\right)+k.
$$
Hence if $k=0$ we get
$$
F\left(B\left(e(A)\right)\right)+F(A)=c_0\Rightarrow
$$
$$
\left(F\circ B\right)\left(e(B(A))\right)+\left(F\circ B\right)(A)=c_0\Leftrightarrow
$$
$$
\left(F\circ B\right)\left(h_0(A)\right)+\left(F\circ B\right)(A)=c_0\Rightarrow
$$
$$
F(B(e(A)))+F(B(e(B(e(A)))))=c_0\Leftrightarrow F(\lambda(A))+F(\lambda(\lambda(A)))=c_0\Leftrightarrow
$$
$$
F(\lambda(A))+F(A)=c_0
$$
For (105) we have
$$
F(B(A))=\frac{w(e(B(A)))}{2\pi i}+c_1=\frac{w(h_0(A))}{2\pi i}+c_1\Rightarrow
$$
$$
F(B(A))=\frac{2\pi i (c_0-2c_1)-w(A)}{2\pi i}+c_1\Leftrightarrow 
$$
$$
2\pi i F(B(A))+w(A)=2\pi i c_1+2\pi i (c_0-2c_1)
$$
which give us immediately (105). The proof of other identities are similar and easy.\\
\\

Now set
$$
c'=c_1-c_0+\frac{c}{2\pi i}.
$$
Hence easily
$$
c'+F\left(B\left(q_A\right)\right)=\frac{c}{2\pi i}-\frac{w(q_A)}{2\pi i}=h_i(A)\Rightarrow
$$
\begin{equation}
c'+F\left(\lambda(A)\right)=h_i(A).
\end{equation}
Hence $h_i(A+1)=h_i(A)$. Also
$$
\int^{\lambda(A)}_{Y(i\infty)}\frac{dt}{P(t)}=\int^{i\infty}_{Y(i\infty)}\frac{dt}{P(t)}+\int^{\lambda(A)}_{i\infty}\frac{dt}{P(t)}=h_i(A)\Rightarrow
$$
$$
c_0-c_1-\int^{A}_{i\infty}\frac{dt}{P(t)}=h_i(A)\Rightarrow
h_i(A)=c''-\int^{A}_{i\infty}\frac{dt}{P(t)}\Rightarrow 
$$
$$
A=c_0-c_1-\int^{h(A)}_{i\infty}\frac{dt}{P(t)}\Rightarrow
h'(A)=-P(h(A)),
$$
where $c''=c_0-c_1$. Also from (27):
$$
w(q_A)=-2\pi i h_i(A)+c
$$
and from (94):
$$
w(e(A))+w\left(h_0(e(A))\right)=c_0-2c_1\Rightarrow
$$
$$
w(e(A))+w\left(e\left(\lambda(A)\right)\right)=c_0-2c_1\Rightarrow
$$
$$
-2\pi i h_i(A)-2\pi i h_i\left(\lambda(A)\right)+2c=c_0-2c_1.
$$
Hence we get the next\\
\\
\textbf{Theorem 10.}\\
We have
\begin{equation}
h_i(A)=c'+F(\lambda(A))\textrm{, }h_i(A+1)=h_i(A),
\end{equation}
\begin{equation}
h'(A)=-P(h(A))\textrm{, }\lambda(A)=h(F(A)),
\end{equation}
\begin{equation}
P(A)h_i'(A)=-1
\end{equation}
and
\begin{equation}
h_i(A)+h_i\left(\lambda(A)\right)=\frac{-c_0+2c_1+2c}{2\pi i}.
\end{equation}
\textbf{Remark.} $h_i(A)$ denotes inversion i.e. $h_i(A)=h^{(-1)}(A)$, $f_i'(A)={f^{(-1)}}{'}(A), \ldots$ etc.\\
\\
\textbf{Theorem 11.}\\
There exists constants $c,c_1'$ such that
\begin{equation}
h_i(A)=\frac{c}{2\pi i}-\frac{w(q_A)}{2\pi i}=\frac{c_1'}{2\pi i}-F(A).
\end{equation}
\\

About the ''shape'' of function $G$, we assume first that $G(F_1(A))$ is analytic and set
$$
H(z):=G\left(F_1\left(z\right)\right).
$$ 
Then
$$
H(z+z_0)=\sum^{\infty}_{k=0}\frac{H^{(k)}(z_0)}{k!}z^k.
$$
Hence we can write
$$
-P(A)=-\frac{1}{q_Aw'(q_A)}=H\left(-\frac{w(q_A)}{2\pi i}+\frac{c}{2\pi i}\right)=
$$
$$
=\sum^{\infty}_{k=0}\frac{H^{(k)}\left(\frac{c}{2\pi i}\right)}{k!}\frac{(-1)^k}{(2\pi i)^k}w(q_A)^k\Rightarrow
$$
$$
-\frac{1}{Aw'(A)}=\sum^{\infty}_{k=0}\frac{H^{(k)}\left(\frac{c}{2\pi i}\right)}{k!}\frac{(-1)^k}{(2\pi i)^k}w(A)^k\Rightarrow
$$
$$
-\frac{{w^{(-1)}}{'}(A)}{w^{(-1)}(A)}=\sum^{\infty}_{k=0}\frac{H^{(k)}\left(\frac{c}{2\pi i}\right)}{k!}\frac{(-1)^k}{(2\pi i)^k}A^k\Rightarrow
$$
$$
-\log w^{(-1)}(A)=c+\sum^{\infty}_{k=0}\frac{H^{(k)}\left(\frac{c}{2\pi i}\right)}{k!}\frac{(-1)^k}{(2\pi i)^k}\frac{A^{k+1}}{k+1}\Rightarrow
$$
$$
w^{(-1)}(A)=\exp\left(-c-A\sum^{\infty}_{k=0}\frac{H^{(k)}\left(\frac{c}{2\pi i}\right)}{(k+1)!}\frac{(-1)^k}{(2\pi i)^k}A^{k}\right)\Rightarrow
$$
$$
f(A)=A\exp\left(c+A\sum^{\infty}_{k=0}\frac{H^{(k)}\left(\frac{c}{2\pi i}\right)}{(k+1)!}\frac{(-1)^k}{(2\pi i)^k}A^{k}\right).
$$
Hence given a function $G$, we can find $f$ setting $H(z)=G(F_1(z))$ and 
\begin{equation}
f(A)=A\exp\left(c+A\sum^{\infty}_{k=0}\frac{H^{(k)}\left(\frac{c}{2\pi i}\right)}{(k+1)!}\frac{(-1)^k}{(2\pi i)^k}A^{k}\right).
\end{equation}
However we have assumed that $f(0)$ is not zero arround 0 ''say'' in $D\subset\textbf{C}$. Hence $G(F_1(A))$ must have a pole. We can write 
$$
G(F_1(A))=-\frac{1}{c-2\pi i A}+P_0(c-2\pi i A),\eqno{(111.1)}
$$
for some function $P_0(A)$. Then
$$
G(F_1(A))=\frac{-1}{c-2\pi i A}+\sum^{\infty}_{k=0}\frac{P_0^{(k)}\left(c\right)(-1)^k (2\pi i)^k}{k!}A^k.
$$
That is because
$$
G\left(F_1\left(\frac{c}{2\pi i}-\frac{A}{2\pi i}\right)\right)=-\frac{1}{A}+\sum^{\infty}_{k=0}\frac{P_0^{(k)}\left(c\right)(-1)^k(2\pi i)^k}{k!}\left(\frac{c}{2\pi i}-\frac{A}{2\pi i}\right)^k=
$$
$$
=-\frac{1}{A}+\sum^{\infty}_{k=0}\frac{P^{(k)}_0\left(c\right)}{k!}(A-c)^k=-\frac{1}{A}+P_0(A),
$$
which is true. Also
$$
-P(A)=-\frac{1}{w(q_A)}+P_0(w(q_A))\Leftrightarrow -\frac{1}{q_Aw'(q_A)}=-\frac{1}{w(q_A)}+P_0(w(q_A))\Leftrightarrow
$$
$$
-\frac{1}{Aw'(A)}=-\frac{1}{w(A)}+P_0(w(A))\Leftrightarrow-\frac{{w^{(-1)}}{'}(A)}{w^{(-1)}(A)}=-\frac{1}{A}+P_0(A)\Leftrightarrow
$$
$$
-\log w^{(-1)}(A)=-\log A+\int^{A}_{C_1}P_0(t)dt+C_0\Leftrightarrow
$$
$$
 \frac{1}{w^{(-1)}(A)}=\frac{e^{C_0}}{A}\exp\left(\int^{A}_{C_1}P_0(t)dt\right)\Leftrightarrow
w^{(-1)}(A)=Ae^{-C_0}\exp\left(-\int^{A}_{C_1}P_0(t)dt\right).
$$
Hence
$$
f(A)=\exp\left(C_0+\int^{A}_{C_1}P_0(t)dt\right)\Leftrightarrow \frac{f'(A)}{f(A)}=P_0(A).
$$
Hence if $f$ is analytic and not zero arround $0$, then so $P_0(A)$ is also analytic and we have the next\\
\\
\textbf{Theorem 12.}\\
Assuming the problem (1),(2),(7), be well defined and in accordance with Ramanujan-Jacobi integral (42), the function  $G(F_1(A))$ must be meromorphic with a single simple pole at $A_0=\frac{c}{2\pi i}$. The constant $c$ is given by $c=w(e(X(0)))$. Moreover it holds
\begin{equation}
G(F_1(A))=\frac{-1}{c-2\pi i A}+\sum^{\infty}_{k=0}\frac{P_0^{(k)}\left(c\right)(-1)^k (2\pi i)^k}{k!}A^k,
\end{equation}
where
\begin{equation}
f(A)=\exp\left(C_0+\int^{A}_{C_1}P_0(t)dt\right)\Leftrightarrow \frac{f'(A)}{f(A)}=P_0(A).
\end{equation}
Hence setting $P_0(A)=\frac{f'(A)}{f(A)}$, then $G$ is given by (111.1) and the oposite.\\
\\

If $D_1=b_1^2-4a_1c_1$ and 
$$
U(x)=U(a_1,b_1;m;x)=
$$
\begin{equation}
=(-1)^{m+1}a_1^{m-1}D_1^{-m+1/2}B_0\left(\frac{-b_1+\sqrt{D_1}-2a_1x}{2\sqrt{D_1}};1-m,1-m\right).
\end{equation} 
Then from [2]:
$$
\exp\left(2\pi i\int^{\omega_2}_{\omega_1}\frac{f_1(t)}{(a_1t^2+b_1t+c_1)^m}dt\right)=\frac{\exp\left(2\pi i \left(h(U(\omega_2)\right)\right)}{\exp\left(2\pi i \left(h(U(\omega_1)\right)\right)}=
$$
$$
=\frac{w^{(-1)}\left(-2\pi i U(\omega_2)+c\right)}{w^{(-1)}\left(-2\pi i U(\omega_1)+c\right)}=
$$
$$
=\frac{(-2\pi i U\left(\omega_2\right)+c)\exp\left(-\int^{-2\pi i U(\omega_2)+c}_{C_1}P_0(t)dt\right)}{(-2\pi i U\left(\omega_1\right)+c)\exp\left(-\int^{-2\pi i U(\omega_1)+c}_{C_1}P_0(t)dt\right)},
$$
where
$$
f_1(A)=-\frac{1}{c-2\pi i U(A)}+P_0\left(c-2\pi i U(A)\right).
$$
But $h'(A)=G(F_1(A))=f_1(U^{(-1)}(A))$. Hence $P_0(A)=\frac{f'(A)}{f(A)}$ and $\frac{w(q)}{f(w(q))}=q$. Therefore
$$
f_1\left(U^{(-1)}\left(\frac{c}{2\pi i}-\frac{A}{2\pi i}\right)\right)=-\frac{1}{A}+P_0(A).
$$
Hence 
$$
f_1\left(U^{(-1)}(A)\right)+\frac{1}{c-2\pi i A}=\sum^{\infty}_{k=0}\frac{P_0^{(k)}\left(c\right)(-1)^k (2\pi i)^k}{k!}A^k\Rightarrow
$$
$$
h(A)=\frac{1}{2\pi i}\log\left(c-2\pi i A\right)+A\sum^{\infty}_{k=0}\frac{P_0^{(k)}(c)(-1)^k(2\pi i)^k}{(k+1)!}A^k+C_1
$$
and
$$
h^{(-1)}(A)=-\frac{w(q_A)}{2\pi i}+\frac{c}{2\pi i}.
$$
Hence we have the next theorem\\
\\
\textbf{Theorem 13.}\\
\textbf{i)} Assume that the function $f_1$ is known and of the form
\begin{equation}
f_1\left(U^{(-1)}\left(\frac{c-A}{2\pi i}\right)\right)=-\frac{1}{A}+P_0\left(A\right),
\end{equation}
where $P_0(A)$ analytic arround 0. Then
$$
\exp\left(2\pi i\int^{\omega_2}_{\omega_1}\frac{f_1(t)}{(a_1t^2+b_1t+c_1)^m}dt\right)=
$$
\begin{equation}
=\frac{(-2\pi i U\left(\omega_2\right)+c)\exp\left(-\int^{-2\pi i U(\omega_2)+c}_{C_1}P_0(t)dt\right)}{(-2\pi i U\left(\omega_1\right)+c)\exp\left(-\int^{-2\pi i U(\omega_1)+c}_{C_1}P_0(t)dt\right)}.
\end{equation}
\textbf{ii)}
\begin{equation}
h(A)=\frac{1}{2\pi i}\log\left(c-2\pi i A\right)+A\sum^{\infty}_{k=0}\frac{P_0^{(k)}(c)(-1)^k(2\pi i)^k}{(k+1)!}A^k+C_1
\end{equation}
and $C_1$ being a constant.\\
\\
\textbf{Theorem 13.1}\\
Given the functions $f_1(A)$ and $U(A)$ ($U$ being that of (115)), such that $f_1\left(U^{(-1)}\left(\frac{c-A}{2\pi i}\right)\right)$ is meromorphic with only simple pole at $A_0=0$ and residue $-1$ i.e. with Laurent expansion (116) and  $P_0(A)=\frac{f'(A)}{f(A)}$ is analytic. Then with the notation of the remarks below we have
\begin{equation}
\int^{-\rho_1-\frac{\sqrt{D_1}}{a_1}\beta(z_2)}_{-\rho_1-\frac{\sqrt{D_1}}{a_1}\beta(z_1)}\frac{f_1(t)}{(a_1t^2+b_1t+c_1)^m}dt
=\frac{1}{2\pi i}\left[\log\left(\frac{c-2\pi i t}{f\left(c-2\pi i t\right)}\right)\right]^{t=\Omega(z_2)}_{t=\Omega(z_1)}.
\end{equation}
The function $\Omega(z)$ is
\begin{equation}
\Omega(z):=(-1)^{m+1}a_1^{m-1}D_1^{(-m+1/2)}\frac{\Gamma(1-m)^2}{\Gamma(2(1-m))}\frac{1}{1-z^2}
\end{equation}
and $\rho_1=\frac{b_1-\sqrt{D_1}}{2a_1}$. The function $f$ can evaluated as in Theorem 6.2 with the help of Theorem 6.1.\\
\textbf{Remarks.}\\
\textbf{i)} Setting
\begin{equation}
B_{\alpha}(z)=\sqrt{B_0\left(z;\alpha,\alpha\right)}\textrm{, }0<\alpha<1,
\end{equation}
the equation
\begin{equation}
i\frac{B_{1-m}\left(1-t\right)}{B_{1-m}\left(t\right)}=z\textrm{, }0<m<1\textrm{, }Im(z)>0,
\end{equation}
have solution $t=\beta(z)$. For this solution holds
\begin{equation}
B_{0}(\beta(z);1-m,1-m)=\frac{\Gamma(1-m)^2}{\Gamma(2(1-m))}\frac{1}{1-z^2}.
\end{equation}
Also
\begin{equation}
U\left(-\rho_1-\frac{\sqrt{D_1}}{a_1}\beta(z)\right)=(-1)^{m+1}a_1^{m-1}D_1^{-m+1/2}\frac{\Gamma(1-m)^2}{\Gamma(2(1-m))}\frac{1}{1-z^2}.
\end{equation}
\textbf{ii)} Hence the beta functions $B_{1-m}(z)$ form (122). For $m$ rational in $0<m<1$, we have numerical evidences that $\beta(z)$ are algebraic numbers when $z=x+i \sqrt{y}$, $x,y$ rationals, with $y>0$. Hence it is of interest to examine these functions. Also it is of interest to reduce the evaluation of general integrals such the left side of (119) with these simple functions.\\
\textbf{iii)} Therorem 13.1 tell us that if $f_1$ is a function such that $f_1\circ U_i$ have simple Laurent expansion, then we can evaluate integral (119) using the analytic part $P_0$ of $f_1\circ U_i$. The evaluation requires only the knowledge of $f$ and $\frac{f'(A)}{f(A)}=P_0(A)$.\\
\textbf{iv)} The problem also related with Ramanujan-Jacobi integrals (see relation (29) and [2]) and holds  $h'(A)=G(F_1(A))=f_1(U_i(A))$. This last equation and Theorem 2 give rise to Lagrange inversion formula, since it holds
$$
h\left(\frac{c}{2\pi i}-\frac{w(q_A)}{2\pi i}\right)=A.   
$$ 
\textbf{v)} By Theorem 6.2.1, the integral (119) can be evaluated if we know the function $y(A)$. Hence asuming the problem (62) of Theorem 6.1 solved we can give evaluations of (119). However, the best solutions of problem (62) until now are equations (66),(66.6.3) and Theorem's 25.5, 25.6, 25.7 below. Hence: knowing $f_1$ we can find $G$. Then from Theorem 6.1, 6.2, we get the inverse of $y(A)$ (relation (66),(66.6.3)). Also using Theorem's 25.5, 25.6, 25.7, we get $y(A)$. From (45) we find $P(A)$. Using (291) of Theorem 25.3 below we find $f(A)$. All these are evaluations not involving integrals and series. We use only simple inversions and one can say that a problem has ''good reduction'' iff it can be described by solving algebraic equations (i.e. inverses of functions).\\ 
\textbf{v)} If we know $G(A)$, then from (67) we get $P_0(A)$ and from (66.6.3) we get $y(A)$. Then from (45) we get $P(A)$ and from (55) we get $w(q)$. Hence from (1) we get $f(A)$ and the integral (119) follows.\\ 
\\
\textbf{Proof.}\\
Assume the Lagrange equation
$$
\frac{w(A)}{f(w(A))}=A.
$$
We find $P_0(A)$ from 
$$
P_0(A)=\frac{f'(A)}{f(A)}.
$$
Then holds the following integral
\begin{equation}
\int^{U_i(h_i(z_2))}_{U_i(h_i(z_1))}\frac{f_1(t)}{(a_1t^2+b_1t+c_1)^m}dt=z_2-z_1
\end{equation}
where
$$
f_1(A)=-\frac{1}{c-2\pi i U(A)}+P_0\left(c-2\pi i U(A)\right).
$$
and
\begin{equation}
w(e(A))=-2\pi ih_i(A)+c.
\end{equation}
Hence
$$
\int^{z_2}_{z_1}\frac{f_1(t)}{(a_1t^2+b_1t+c_1)^m}dt=\frac{1}{2\pi i}\log\left(\frac{c-2\pi i U(z_2)}{c-2\pi i U(z_1)}\right)+
$$
$$
+\sum^{\infty}_{k=0}\frac{P_0^{(k)}(c)(-2\pi i)^k}{(k+1)!}\left(U(z_2)^{k+1}-U(z_1)^{k+1}\right).\eqno{(126.1)}
$$
Hence if we set in $A_2,A_1$ the values
$$
A_2=-\rho_1-\frac{\sqrt{D_1}}{a_1}\beta(z_2)\textrm{, }A_1=-\rho_1-\frac{\sqrt{D_1}}{a_1}\beta(z_1),
$$
we get 
\scriptsize 
$$
\int^{-\rho_1-\frac{\sqrt{D_1}}{a_1}\beta(z_2)}_{-\rho_1-\frac{\sqrt{D_1}}{a_1}\beta(z_1)}\frac{f_1(t)}{(a_1t^2+b_1t+c_1)^m}dt=\frac{1}{2\pi i}\log\left(\frac{c-2\pi i U(A_2)}{c-2\pi i U(A_1)}\right)+
$$
$$
+\sum^{\infty}_{k=0}\frac{P_0^{(k)}(c)(-2\pi i)^k}{(k+1)!}(-1)^{(m+1)(k+1)}a_1^{(m-1)(k+1)}D_1^{(-m+1/2)(k+1)}\frac{\Gamma(1-m)^{2(k+1)}}{\Gamma(2(1-m))^{k+1}}\left[\frac{1}{(1-t^2)^{k+1}}\right]^{t=z_2}_{t=z_1}.
$$
\normalsize
Hence if $P^{(k)}_{0}(c)=P^{(k+1)}_1(c)$, then
$$
\int^{-\rho_1-\frac{\sqrt{D_1}}{a_1}\beta(z_2)}_{-\rho_1-\frac{\sqrt{D_1}}{a_1}\beta(z_1)}\frac{f_1(t)}{(a_1t^2+b_1t+c_1)^m}dt=\frac{1}{2\pi i}\log\left(\frac{c-2\pi i U(A_2)}{c-2\pi i U(A_1)}\right)-
$$
$$
-\frac{1}{2\pi i}\left[P_1\left(c-2\pi i(-1)^{m+1}a_1^{m-1}D_1^{(-m+1/2)}\frac{\Gamma(1-m)^2}{\Gamma(2(1-m))}\frac{1}{1-t^2}\right)\right]^{t=z_2}_{t=z_1}\Rightarrow
$$
$$
\int^{-\rho_1-\frac{\sqrt{D_1}}{a_1}\beta(z_2)}_{-\rho_1-\frac{\sqrt{D_1}}{a_1}\beta(z_1)}\frac{f_1(t)}{(a_1t^2+b_1t+c_1)^m}dt=\frac{1}{2\pi i}\log\left(\frac{c-2\pi i U(A_2)}{c-2\pi i U(A_1)}\right)-
$$
$$
-\frac{1}{2\pi i}\int^{-2\pi i(-1)^{m+1}a_1^{m-1}D_1^{(-m+1/2)}\frac{\Gamma(1-m)^2}{\Gamma(2(1-m))}\frac{1}{1-z_2^2}}_{-2\pi i(-1)^{m+1}a_1^{m-1}D_1^{(-m+1/2)}\frac{\Gamma(1-m)^2}{\Gamma(2(1-m))}\frac{1}{1-z_1^2}}P_0(t+c)dt\Rightarrow
$$
$$
\int^{-\rho_1-\frac{\sqrt{D_1}}{a_1}\beta(z_2)}_{-\rho_1-\frac{\sqrt{D_1}}{a_1}\beta(z_1)}\frac{f_1(t)}{(a_1t^2+b_1t+c_1)^m}dt=\frac{1}{2\pi i}\log\left(\frac{c-2\pi i U(A_2)}{c-2\pi i U(A_1)}\right)-
$$
$$
-\frac{1}{2\pi i}\int^{-2\pi i(-1)^{m+1}a_1^{m-1}D_1^{(-m+1/2)}\frac{\Gamma(1-m)^2}{\Gamma(2(1-m))}\frac{1}{1-z_2^2}}_{-2\pi i(-1)^{m+1}a_1^{m-1}D_1^{(-m+1/2)}\frac{\Gamma(1-m)^2}{\Gamma(2(1-m))}\frac{1}{1-z_1^2}}P_0(t+c)dt=
$$
$$
=\frac{1}{2\pi i}\log\left(\frac{c-2\pi i U(A_2)}{c-2\pi i U(A_1)}\right)-
$$
$$
-\frac{1}{2\pi i}\log\left(\frac{f\left(c-2\pi i(-1)^{m+1}a_1^{m-1}D_1^{(-m+1/2)}\frac{\Gamma(1-m)^2}{\Gamma(2(1-m))}\frac{1}{1-z_2^2}\right)}{f\left(c-2\pi i(-1)^{m+1}a_1^{m-1}D_1^{(-m+1/2)}\frac{\Gamma(1-m)^2}{\Gamma(2(1-m))}\frac{1}{1-z_1^2}\right)}\right).
$$
\\
\textbf{Example.}\\
If $P_0(A)=1$, then we have $f(A)=C_0e^A$ and $w_i(A)=C_0^{-1}Ae^{-A}$, with $w(A)=-W\left(-C_0 A\right)$ and $W(A)$ is the Lambert's function, $C_0=e^{-C}$.
$$
f_1(U_i(A))=-\frac{1}{c-2\pi i A}+1\Leftrightarrow f_1(A)=-\frac{1}{c-2\pi i U(A)}+1.
$$
Also then
$$
\exp\left(2\pi i \int^{\omega}_{0}\frac{f_1(t)}{(a_1t^2+b_1t+c_1)^m}dt\right)=e(U(\omega))\left(1-2\pi i c^{-1} U(\omega)\right).
$$
On the other hand we have $G(F_1(A))=-\frac{1}{c-2\pi i A}+1$. Hence
$$
G(A)=-\frac{1}{c-2\pi i F_1^{(-1)}(A)}+1.
$$
The function $h_0(A)$ is such that $w(A)+w(h_0(A))=c_{11}$. Hence
$$
h_0(A):=e^{C-c_{11}-W(-A e^{-C})}\left(c_{11}+W(-A e^{-C})\right)\textrm{, }c_{11}=2\pi i (c_0-c_1)
$$
and indeed holds $h_0(h_0(A))=A$, ($W(x)$ is the Lambert's $W$ function). The function $y(A)$ is
$$
y(A)=F_1\left(\frac{W\left(-qe^{-C}\right)}{2\pi i}+\frac{c}{2\pi i}\right)\textrm{, }q=e(A).
$$
$$
F(A)=\frac{c_1-c}{2\pi i}-\frac{W\left(qe^{-C}\right)}{2\pi i}\textrm{, }q=e(A)
$$
$$
P(A)=-\frac{1+W(-qe^{-C})}{W\left(-qe^{-C}\right)}\textrm{, }q=e(A).
$$
Also from (118)
$$
h(A)=\frac{1}{2\pi i}\log\left(A-\frac{c}{2\pi i}\right)+A+C_1
$$
and from (24),(12) 
$$
Y(z)=\frac{1}{2\pi i}\log\left(\sqrt[3]{2}B_0\left(m^*(2z)^2;\frac{1}{6},\frac{2}{3}\right)-\frac{c}{2\pi i}\right)+\sqrt[3]{2}B_0\left(m^*(2z)^2;\frac{1}{6},\frac{2}{3}\right)+C_1.
$$
\\
\textbf{Note.} Solving equation (9) with ''Mathematica'' program (I have ''Mathematica 11'') requires some extra care when using the constants. Also mathematica does not recognizes $W(xe^x)=x$ and it is beter to use $6A^{1/3} {}_2F_{1}\left(\frac{1}{6},\frac{1}{3};\frac{7}{6};A^2\right)$ in place of $B_0\left(A^2;\frac{1}{6},\frac{2}{3}\right)$. A beter example is to take $P_0(t)=\frac{1}{1+t}$, which is equivalent to $f(A)=C(A+1)$.

\section{The real analog}

Going from the complex to the real analog (see [3]) here we have an equation
\begin{equation}
\frac{w(q)}{f(w(q))}=q,
\end{equation}
with $f(A)$ analytic and $f(0)\neq 0$ arround 0. The equation (127) have solution
\begin{equation}
w(q)=\sum^{\infty}_{n=1}c_nq^n\textrm{, }q=e^{-\pi\sqrt{A}}\textrm{, }A>0.
\end{equation}
Then if
\begin{equation}
a_n=c_n n,
\end{equation}
we will study all the functions in which
\begin{equation}
\frac{1}{P(A)}=\sum^{\infty}_{n=1}a_nq^n\textrm{, }q=e^{-\pi\sqrt{A}}\textrm{, }A>0,
\end{equation}
in the sence that $P(A)$ defines a function $X(A)$ such that X(A) is solution of the equation (137) below and is connected with inversion problem in [3] and [2]. Hence due to the connection (129) the class of all functions $P(A)$ is very wide. A first result is
\begin{equation}
\int \frac{1}{qP(A)}dq=w(q)+c 
\end{equation}
and
\begin{equation}
P(A)=\frac{1}{qw'(q)}\textrm{, }q=e^{-\pi \sqrt{A}}\textrm{, }A>0.
\end{equation}
Also
\begin{equation}
-2P(A)h_i'(A)=1
\end{equation}
and
$$
h_i(A)=c+\frac{1}{\pi^2}\sum^{\infty}_{n=1}\frac{a_n}{n^2}q^n+\frac{\sqrt{A}}{\pi}\sum^{\infty}_{n=1}\frac{a_n}{n}q^n=
$$
$$
=c+\pi^{-2} \int q^{-1}w(q)dq-\pi^{-2}w(q)\log q=
$$
$$
=c+\pi^{-2} w(q)\log q-\pi^{-2}\int w'(q)\log(q)dq-\pi^{-2}w(q)\log q\Rightarrow
$$
\begin{equation}
h_i(A)=c-\pi^{-2}\int w'(q)\log(q) dq
\end{equation}
and
\begin{equation}
h_i'(A)=-\frac{1}{2}e^{-\pi\sqrt{A}}w'\left(e^{-\pi\sqrt{A}}\right).
\end{equation}
The function $X(A)$ is given from
\begin{equation}
X(A)=h\left(\frac{1}{\sqrt[3]{4}}B_0\left(A^2;\frac{1}{6},\frac{2}{3}\right)\right)
\end{equation}
and satisfies the equation
\begin{equation}
X'(A)+\frac{2^{4/3}}{A^{2/3}(1-A^2)^{1/3}}P(X(A))=0,
\end{equation}
which is equivalent to
\begin{equation}
h'(A)+2P(h(A))=0.
\end{equation}
The function
\begin{equation}
Y(r)=X(k_r),
\end{equation}
satisfies
\begin{equation}
Y'(r)-\frac{\pi}{\sqrt{r}}\eta\left(i\sqrt{r}/2\right)^4P(Y(r))=0.
\end{equation}
Also if
\begin{equation}
F(A)=\int^{A}_{X(0)}\frac{dt}{P(t)},
\end{equation}
then
\begin{equation}
F(Y(r))=-\frac{2}{\sqrt[3]{4}}B_0\left(k_r^2;\frac{1}{6},\frac{2}{3}\right).
\end{equation}
and
\begin{equation}
F(Y(4r))+F\left(Y\left(\frac{4}{r}\right)\right)=c_0=-\frac{\sqrt{3}\Gamma\left(\frac{1}{3}\right)}{\pi\sqrt[3]{2}}
\end{equation}
Moreover if $G(A)$ is the function related with the Ramanujan-Jacobi inversion problem, then setting
\begin{equation}
h_1(t):=\left(\frac{1}{h_i'(.)}\right)^{(-1)}(t),
\end{equation}
we have
\begin{equation}
G_i(x)=F_1\left(\int^{x}_{c}\frac{h_1'(t)}{t}dt\right)
\end{equation}
and
\begin{equation}
G^{(-1)}(-2x)=F_1\left(\frac{-1}{2}\int^{P_i(x)}_{c}\frac{dt}{P(t)}\right).
\end{equation}
\\
\textbf{Theorem 14.}
\begin{equation}
w'\left(q\right)q=-2h_i'(A)\textrm{, }q=e^{-\pi\sqrt{A}}\textrm{, }A>0.
\end{equation}
Or equivalent
\begin{equation}
h_i(A)=c-\frac{1}{2}\int w'(q)qdA=c'-\pi^{-2}\int w'(q)\log(q)dq=c'-\pi^{-2}C(q).
\end{equation}
\\

If $q=e^{-\pi\sqrt{A}}$, then
$$
5\int^{y(A)}_{0}\frac{dt}{t\sqrt[6]{t^{-5}-11-t^5}}=c'-\pi^{-2}\int w'(q)\log(q) dq=
$$
$$
=c-\frac{1}{2}\int w'(q)qdA\Rightarrow
$$
\begin{equation}
F_1\left(c-\frac{1}{2}\int w'(q)q dA\right)=y(A),
\end{equation}
where
$$
5\int^{F_1(A)}_{0}\frac{dt}{t\sqrt[6]{t^{-5}-11-t^5}}=A.
$$
\textbf{Note.} In view of [2] we have for $A$ real and positive
$$
F_1(A)=R\left(e^{-\pi\sqrt{m_0(A)}}\right),\eqno{(149.1)}
$$
where $m_0(A)$ is the inverse function of
$$
2^{-2/3}B_0\left(k_A^2;\frac{1}{6},\frac{2}{3}\right)
$$
and $R(q)$ is the Rogers-Ramanujan continued fraction. The function $k_r$ being the elliptic singular modulus i.e: 
$$
k_r=\left(\frac{\theta_2(q)}{\theta_3(q)}\right)^2\textrm{, }q=e^{-\pi\sqrt{r}}\textrm{, }r>0.
$$
But $G(y(A))=1/h_i'(A)$. Hence
\begin{equation}
G\left(F_1\left(c-\frac{1}{2}\int w'(q)qdA\right)\right)=-2P(A).
\end{equation}
\\

Assume (in the same way as we did above) that 
$$
G\left(F_1\left(A\right)\right)=H(A),
$$
where $H(A)$ is analytic. Also set 
\begin{equation}
C(q)=\int w'(q)\log(q) dq=\int^{q}_{}w'(t)\log(t) dt
\end{equation}
and $f(A)=e^{-P_1'(A)}$ analytic and not zero at the origin. Then 
$$
C'(A)=w'(A)\log A\Rightarrow C'(w_i(A))w_i'(A)=w'(w_i(A))\log(w_i(A))w_i'(A)\Rightarrow
$$
$$
C(w_i(A))=\int\log(w_i(A))dA+c_1=\int\log\left(\frac{A}{f(A)}\right)dA+c_1=
$$
\begin{equation}
=\int\log(A)dA+P_1(A)+c_1.
\end{equation}
We have
$$
H(A)=\sum^{\infty}_{k=0}\frac{H^{(k)}(0)}{k!}A^k.
$$
Setting $A\rightarrow c-\frac{1}{2}\int w'(q)qdA$, we have
$$
G\left(F_1\left(c-\frac{1}{2}\int w'(q)qdA\right)\right)=\sum^{\infty}_{k=0}\frac{H^{(k)}(0)}{k!}\left(c-\frac{1}{2}\int w'(q)qdA\right)^k\Rightarrow
$$
$$
-2P(A)=\sum^{\infty}_{k=0}\frac{H^{(k)}(0)}{k!}\left(c'-\pi^{-2}\int w'(q)qdq\right)^k\Rightarrow
$$
$$
\frac{-2}{w'(q)q}=\sum^{\infty}_{k=0}\frac{H^{(k)}(0)}{k!}\left(c'-\pi^{-2}C(q)\right)^k\Rightarrow
$$
$$
-\frac{2}{w'(A)A}=\sum^{\infty}_{k=0}\frac{H^{(k)}(0)}{k!}\left(c'-\pi^{-2}C(A)\right)^{k}\Rightarrow
$$
$$
-2\frac{w_i'(A)}{w_i(A)}=\sum^{\infty}_{k=0}\frac{H^{(k)}(0)}{k!}\left(c'-\pi^{-2}C(w_i(A))\right)^{k}\Rightarrow
$$
$$
\frac{w_i'(A)}{w_i(A)}=-\frac{1}{2}\sum^{\infty}_{k=0}\frac{H^{(k)}(0)}{k!}\left(c'-c_1-\pi^{-2}\int\log(A)dA-\pi^{-2}P_1(A)\right)^{k}.
$$
Now if $\xi$ is positive constant and $x$ positive variable with $0<x<\xi$, then
$$
\log w_i(x)-\log w_i(\xi)=
$$
$$
-\frac{1}{2}\int^{x}_{\xi}\left\{\sum^{\infty}_{k=0}\frac{H^{(k)}(0)}{k!}\left(c'-c_1-\pi^{-2}\int\log(A)dA-\pi^{-2}P_1(A)\right)^{k}\right\}dA\Rightarrow
$$
$$
w_i(x)=w_i(\xi)\times
$$
$$
\times\exp\left(-\frac{1}{2}\int^{x}_{\xi}\left\{\sum^{\infty}_{k=0}\frac{H^{(k)}(0)}{k!}\left(c'-c_1-\pi^{-2}\int \log(A)dA-\pi^{-2}P_1(A)\right)^{k}\right\}dA\right).
$$
Finaly
$$
f(x)=\frac{x f(\xi)}{\xi}\times
$$
$$
\times\exp\left(\frac{1}{2}\int^{x}_{\xi}\left\{\sum^{\infty}_{k=0}\frac{H^{(k)}(0)}{k!}\left(c'-c_1-\pi^{-2}\int \log(A)dA-\pi^{-2}P_1(A)\right)^{k}\right\}dA\right).
$$
Since the integral $\int \log(A)dA$ is continuous and bounded in $[0,\xi]$ and $P_1(A)$ analytic, we have $f(0)=0$, which is imposible. Hence $G\left(F_1(A)\right)$ is not analytic. However if we assume that 
$$
G\left(F_1\left(A\right)\right)+\frac{2}{L(A)}=P^{*}_0(A),
$$
where $L$ is a function such that
$$
L\left(c-\frac{1}{2}\int w'(q)qdA\right)=w(q)
$$
and $P^{*}_0(A)$ analytic. Then
$$
G\left(F_1\left(c-\frac{1}{2}\int w'(q)qdA\right)\right)+\frac{2}{w(q)}=P^{*}_{0}\left(c-\frac{1}{2}\int w'(q)qdA\right)\Rightarrow
$$
$$
-2P(A)+\frac{2}{w(q)}=P_0^{*}\left(c'-\pi^{-2}\int w'(q)\log (q)dq\right)\Rightarrow
$$
$$
-\frac{2}{w'(q)q}+\frac{2}{w(q)}=P^{*}_{0}\left(c'-\pi^{-2}C(q)\right)\Rightarrow
$$
$$
-2\frac{w_i'(A)}{w_i(A)}+2\frac{1}{A}=P^{*}_{0}\left(c'-\pi^{-2}C(w_i(A))\right)\Rightarrow
$$
$$
\left[-2\log(w_i(t))+2\log t\right]^{t=A}_{t=\xi}=\int^{A}_{\xi}P^{*}_{0}\left(c'-\pi^{-2}C(w_i(t))\right)dt\Rightarrow
$$
$$
-2\log\left(\frac{w_i(A)}{A}\right)+2\log\left(\frac{w_i(\xi)}{\xi}\right)=
$$
$$
=\int^{A}_{\xi}P^{*}_{0}\left(c'-c_1-\pi^{-2}\int \log(t)dt-\pi^{-2}P_1(t)\right)dt\Rightarrow
$$
$$
\log(f(A))-\log(f(\xi))=
$$
$$
=\frac{1}{2}\int^{A}_{\xi}P_0^{*}\left(c'-c_1-\pi^{-2}\int \log(t)dt-\pi^{-2}P_1(t)\right)dt\Rightarrow
$$
$$
f(A)=f(\xi)\exp\left(\frac{1}{2}\int^{A}_{\xi}P_0^{*}\left(c'-c_1-\pi^{-2}\int \log(t)dt-\pi^{-2}P_1(t)\right)dt\right).
$$
Hence we get the next\\
\\
\textbf{Theorem 15.}\\
If we assume the problem
\begin{equation}
\frac{w(q)}{f(w(q))}=q\textrm{, }q=e^{-\pi\sqrt{r}}\textrm{, }r>0,
\end{equation}
where $f$ is analytic arround $0$ and $f(0)\neq 0$ and assume the solution is
\begin{equation}
w(q)=\sum^{\infty}_{n=1}c_nq^n.
\end{equation}
Setting $a_n=nc_n$, we define
\begin{equation}
\frac{1}{P(A)}=\sum^{\infty}_{n=1}a_nq^n.
\end{equation}
 Hence, we define a connection of the Lagrange problem with the Hauptmodul and the Ramanujan-Jacobi problem as described in the above notes of present article (see also [2],[3]). The connection of Ramanujan-Jacobi problem and the Hauptmodul general problem is
\begin{equation}
G_i(-2x)=F_1\left(\frac{-1}{2}\int^{P_i(x)}_{c}\frac{dt}{P(t)}\right).
\end{equation}
Then in order the above problem to be well defined, function $G(F_1(A))$ must be of the form 
\begin{equation}
G\left(F_1\left(A\right)\right)=-\frac{2}{L(A)}+P^{*}_{0}(A),
\end{equation}
where $P^{*}_{0}(A)$ is analytic in a interval containnig $0$ and $L(A)$ must satisfies
\begin{equation}
L\left(c-\frac{1}{2}\int w'(q)qdA\right)=w(q).
\end{equation}
Also if $\xi$ is suitable positive constant and $x>0$, then setting $f(A)=e^{-P_1'(A)}$, we have
\begin{equation}
f(A)=f(\xi)\exp\left(\frac{1}{2}\int^{A}_{\xi}P_0^{*}\left(c'-c_1-\pi^{-2}\int \log(t)dt-\pi^{-2}P_1(t)\right)dt\right).
\end{equation}
\\

The function $L(A)$ can be written as
$$
L_i(w(q))=c-\frac{1}{2}\int w'(q)q dA\Rightarrow L_i'(w(q))w'(q)q\frac{-\pi}{2\sqrt{A}}=-\frac{1}{2}w'(q)q\Leftrightarrow
$$
$$
L_i'(w(q))\frac{\pi}{\sqrt{A}}=1\Leftrightarrow -\pi^2 L_i'(w(q))=-\pi\sqrt{A}=\log q\Leftrightarrow
$$
$$
L_i'(q)=-\frac{1}{\pi^2}\log(w_i(q))\Leftrightarrow L_i'(A)=-\frac{1}{\pi^2}\log\left(\frac{A}{f(A)}\right).
$$
Hence we get the next\\
\\
\textbf{Theorem 16.}
\begin{equation}
L_i(A)=-\pi^{-2}\int \log\left(\frac{A}{f(A)}\right)dA+c.
\end{equation}
\\
\textbf{Theorem 17.}\\
If we set 
\begin{equation}
S(A):=\frac{\pi^{-2}}{2}P^{*}_{0}\left(A\right),
\end{equation}
then
\begin{equation}
-\frac{L''(A)}{L'(A)^3}+\frac{\pi^{-2}}{L(A)}=S(A).
\end{equation}
\\
\textbf{Proof.}\\
We have
$$
-\frac{2}{w'(q)q}+\frac{2}{w(q)}=P_0^{*}\left(c'-\pi^{-2}C(q)\right)\Leftrightarrow
$$
$$
-\frac{2}{w'(q)q}+\frac{2}{w(q)}=P_0^{*}\left(L_i(w(q))\right)\Leftrightarrow
$$
$$
-2\log q+2\log(w(q))=\int^{w(q)}_{c_1}P_0^{*}\left(L_i(t)\right)dt\Leftrightarrow
$$
$$
-2\log w_i(A)+2\log A=\int P_0^{*}\left(L_i(A)\right)dA\Leftrightarrow
$$
$$
2\log(f(A))=\int P^{*}_0(L_i(A))dA\Leftrightarrow
$$
$$
2\log(f(A))=\int P^{*}_{0}\left(\pi^{-2}\int \log\left(\frac{f(A)}{A}\right)\right)dA.
$$
Set $S(A):=\frac{\pi^{-2}}{2}P^{*}_{0}\left(A\right)$. Then 
$$
L_i'(A)=\int S(L_i(A))dA-\pi^{-2}\log A\Rightarrow
$$
\begin{equation}
L_i''(A)=S(L_i(A))-\frac{\pi^{-2}}{A}.
\end{equation}
Hence if $u(x):=L_i(x)$, then
$$
\frac{d}{dA}\left(u''(A) A\right)=\frac{d}{dA}(S(u(A))A)\Rightarrow
$$
\begin{equation}
u''(A)-S(u(A))+A\left[u'''(A)-S'(u(A))u'(A)\right]=0.
\end{equation}
A solution of (164) is
$$
u_i(A)=\int^{A}_{c_1}\frac{dt}{\sqrt{c_2+2\int^{t}_{c_3}S(t_1)dt_1}}.
$$
However if we set $A\rightarrow L(A)$ in (163), then
$$
-\frac{L''(A)}{L'(A)^3}+\frac{\pi^{-2}}{L(A)}=S(A).
$$
\textbf{Remarks.}\\
Hence if exists $S_1(A)$  such that
$$
-\frac{u_i''(A)}{u_i'(A)^3}+\frac{1}{u_i(A)}=S_1\left(u_i(A)\right)=S(A),
$$
we get solving the first equality
$$
u_i(A)=\int^{A}_{c_0}\frac{dt}{c_1-\log t+\int S_1(t)dt}.
$$
Hence
$$
S_1\left(\int^{A}_{c_0}\frac{dt}{c_1-\log t+\int S_1(t)dt}\right)=S(A).
$$
Now set the functions $S,S_0$ such that
\begin{equation}
u_1''(A)=S(u_1(A))-\frac{1}{A}
\end{equation}
and 
\begin{equation}
u_2''(A)=S_0(u_2(A)).
\end{equation}
Hence
$$
u_2^{(-1)}(x)=\int^{x}_{c_1}\frac{dt}{\sqrt{2\int^{t}_{c_2}S_0(t_1)dt_1}}.
$$
But equation (165) can be written as  
$$
-\frac{{u_{1}^{(-1)}}{''}(A)}{{u_1^{(-1)}}{'}(A)^3}+\frac{1}{u_1^{(-1)}(A)}=S(A)=S_1(u^{(-1)}_1(A))\eqno{(166.1)}
$$
and have solution
$$
u_1^{(-1)}(A)=\int^{A}_{c_0}\frac{dt}{c_1-\log t+\int S_1(t)dt}.
$$
However $u_1^{(-1)}(A)=u_2^{(-1)}(A)$. Hence
$$
\int^{A}_{c_1}\frac{dt}{\sqrt{2\int^{t}_{c_2}S_0(t_1)dt_1}}=\int^{A}_{c_0}\frac{dt}{c_1-\log t+\int S_1(t)dt}\Rightarrow
$$
$$
S_1\left(\int^{A}_{c_0}\frac{dt}{c_1-\log t+\int S_1(t)dt}\right)=S(A)=S_1\left(\int^{A}_{c_1}\frac{dt}{\sqrt{2\int^{t}_{c_2}S_0(t_1)dt_1}}\right).
$$
The solution of (165) is not an easy problem and it might be unsolved.\\ However we can write
$$
L(t_0)=\int^{t_0}_{c_2}\frac{dt_1}{\sqrt{c_3+2\int^{t_1}_{c_4}(S(t_2)-S_1(t_2))dt_2}},
$$
where
$$
S_1(A)=\frac{C_1}{A-C_2}\int^{A}_{c_1}\frac{dt_0}{L(t_0)}.
$$
Hence
$$
S_1(A)=\frac{C_1}{A-C_2}\int^{A}_{c_1}\frac{dt_0}{\int^{t_0}_{c_2}\frac{dt_1}{\sqrt{c_3+2\int^{t_1}_{c_4}(S(t_2)-S_1(t_2))dt_2}}}.
$$
Setting $S_1=S_1(A)$, we get
$$
S_1=S_1(A)=T(S_1,A):=\frac{C_1}{A-C_2}\int^{A}_{c_1}\frac{dt_0}{\int^{t_0}_{c_2}\frac{dt_1}{\sqrt{c_3+2\int^{t_1}_{c_4}(S(t_2)-S_1(t_2))dt_2}}}.\eqno{(166.2)}
$$
Hence
$$
S_1(A)=\frac{C_1}{A-C_2}\int^{A}_{c_1}\frac{dt_0}{\int^{t_0}_{c_2}\frac{dt_1}{\sqrt{c_3+2\int^{t_1}_{c_4}\left(S(t_2)-\left[\frac{C_1}{t_2-C_2}\int^{t_2}_{c_1}\frac{dt_3}{\int^{t_3}_{c_2}\frac{dt_4}{\sqrt{c_3+2\int^{t_4}_{c_4}(S(t_5)-S_1(t_5))dt_5}}}\right]\right)dt_2}}}.
$$
Consequently
$$
S_1=S_1(A)=T(S_1,A)=T(T(S_1,A),A)=T(T(T(S_1,A),A),A)=\ldots=
$$
$$
=T\circ T\circ T\circ\ldots.\eqno{(166.3)}
$$
Hence we get the next:\\
\\
\textbf{Theorem 17.1}\\
A solution of
$$
-\frac{u''(A)}{u'(A)^3}+\frac{1}{u(A)}=S(A),\eqno{(166.4)}
$$
is
$$
u(A)=L_i(A)=L^{(-1)}(A)=S^{(-1)}(S_1(A)),\eqno{(166.5)}
$$
where
$$
S_1(A)=T\circ T\circ T\circ\ldots.\eqno{(166.6)}
$$
\\

Now assume that $m$ is rational and set
$$
B_{\alpha}(x):=\sqrt{B_0(x,\alpha,\alpha)}=\sqrt{\int^{x}_{0}(t-t^2)^{\alpha-1}dt}.
$$ 
It is known that 
\begin{equation}
\int^{z_2}_{z_1}\frac{dt}{(at^2+bt+c)^m}=U(a,b,c;m;z_2)-U(a,b,c;m;z_1).
\end{equation}
Also if
\begin{equation}
A_1=-\rho_1-\frac{\sqrt{D}}{a}\beta_{r_1}\textrm{, }A_2=-\rho_1-\frac{\sqrt{D}}{a}\beta_{r_2},
\end{equation}
where $\rho_1=\frac{b-\sqrt{D}}{2a}$ and $\beta_{r_1}$, $\beta_{r_2}$ are solutions of
\begin{equation}
\frac{B_{1-m}\left(1-\beta_{r_{1,2}}\right)}{B_{1-m}\left(\beta_{r_{1,2}}\right)}=\sqrt{r_{1,2}},
\end{equation}
(if $r$ is positive rational, then $\beta_r$ is algebraic), we have
$$
\int^{A_2}_{A_1}\frac{dt}{(at^2+bt+c)^m}=(-1)^{m+1}a^{m-1}D^{-m+1/2}B_{1-m}\left(\beta_{r_2}\right)^2-
$$
$$
-(-1)^{m+1}a^{m-1}D^{-m+1/2}B_{1-m}\left(\beta_{r_1}\right)^2.
$$
But one can easily see that 
$$
B_{\alpha}^2(z)+B_{\alpha}(1-z)^2=\int^{1}_{0}(t(1-t))^{\alpha}dt=\frac{\Gamma(\alpha)^2}{\Gamma(2\alpha)}.
$$
Setting $z=\beta_r$ in the above formula we have
$$
\frac{B_{\alpha}(1-\beta_r)}{B_{\alpha}(\beta_r)}=\sqrt{r}
$$
and
\begin{equation}
B_{\alpha}(\beta_r)=\sqrt{\frac{\Gamma(\alpha)^2}{\Gamma(2\alpha)(r+1)}}
\end{equation}
and also
\begin{equation}
B_{\alpha}\left(\beta_{n^2r}\right)=\sqrt{\frac{r+1}{n^2r+1}}B_{\alpha}(\beta_r).
\end{equation}
Hence we get the next\\
\\
\textbf{Theorem 18.}\\
If $r_1,r_2$ are rational and $A_1,A_2$ are that of (168) with $\beta_{r_1},\beta_{r_2}$ the algebraic solutions of (169), then
\begin{equation}
\int^{A_2}_{A_1}\frac{dt}{(at^2+bt+c)^m}=(-1)^{m+1}a^{m-1}D^{-m+1/2}\frac{\Gamma(1-m)^2}{\Gamma(2(1-m))}\left(\frac{1}{r_2+1}-\frac{1}{r_1+1}\right).
\end{equation}
Assuming $r_2=r$ and $r_1=+\infty$ we have
\begin{equation}
\int^{-\rho_1-\frac{\sqrt{D}}{a}\beta_r}_{-\rho_1}\frac{dt}{(at^2+bt+c)^m}=(-1)^{m+1}a^{m-1}D^{-m+1/2}\frac{\Gamma(1-m)^2}{\Gamma(2(1-m))}\frac{1}{r+1}.
\end{equation}
\\

But it holds (see [2]) $h'(A)=G(F_1(A))=f_1(U_i(A))$. Hence from Theorem 16 we have
$$
h'(A)=-\frac{2}{L(A)}+P_0^{*}(A)\Rightarrow h(U(A))=-2\int^{U(A)}_{c}\frac{dt}{L(t)}+\int^{U(A)}_{c}P^{*}_{0}(t)dt,
$$
where $U(A)$ is that of (116). Hence assuming $P_0^{*}(A)$ is given analytic function, if
\begin{equation}
f_1(t)=-\frac{2}{L(U(t))}+P_0^{*}(U(t)),
\end{equation}
we have
\begin{equation}
\int^{A_2}_{A_1}\frac{f_1(t)}{(a_1t^2+b_1t+c_1)^m}dt=-2\int^{U(A_2)}_{U(A_1)}\frac{dt}{L(t)}+\int^{U(A_2)}_{U(A_1)}P^{*}_{0}(t)dt,
\end{equation}
where $A_1,A_2$ may be arbitrary. The function $L(A)$ is determined from Theorem 16 equation (160) and $P^{*}_{0}(A)$ from Theorem 17. Hence for $f_1(A)$ we can evaluate the integral (175). Hence as a special case
$$
\int^{-\rho_1-\frac{\sqrt{D_1}}{a_1}\beta_{r_2}}_{-\rho_1-\frac{\sqrt{D_1}}{a_1}\beta_{r_1}}\frac{f_1(t)}{(a_1t^2+b_1t+c_1)^m}dt=-2\int^{U(A_2)}_{U(A_1)}\frac{dt}{L(t)}+
$$
\begin{equation}
+\int^{U(A_2)}_{U(A_1)}P^{*}_{0}(t)dt,
\end{equation}
where
\begin{equation}
U(A_{1,2})=(-1)^{m+1}a_1^{m-1}D_1^{-m+1/2}\frac{\Gamma(1-m)^2}{\Gamma(2(1-m))}\frac{1}{r_{1,2}+1}.
\end{equation}
Continuing from Theorem 17 we can write
$$
\int^{A_2}_{A_1}\frac{f_1(t)}{(a_1t^2+b_1t+c_1)^m}dt=
$$
$$
=-2\int^{U(A_2)}_{U(A_1)}\frac{dt}{L(t)}+\int^{U(A_2)}_{U(A_1)}2\pi^{2}\left(-\frac{L''(t)}{L'(t)^3}+\frac{\pi^{-2}}{L(t)}\right)dt=
$$
$$
=-2\pi^2\int^{U(A_2)}_{U(A_1)}\frac{L''(t)}{L'(t)^3}dt=\pi^2\left(\frac{1}{L'(U(A_2))^2}-\frac{1}{L'(U(A_1))^2}\right).
$$
Also
$$
f_1(A)=\pi^2\left(\frac{d}{dt}\frac{1}{L'(t)^2}\right)_{t=U(A)}.
$$
The above will help us to prove the next:\\
\\
\textbf{Theorem 19.}\\
Assume that $f_1$ is any smooth function of the form
\begin{equation}
f_1(A)=\pi^2\left(\frac{d}{dt}\frac{1}{L'(t)^2}\right)_{t=U(A)}.
\end{equation}
Knowing $L(A)$, we can assume that $R_{1,2}$ are solutions of the equation
\begin{equation}
c-\frac{1}{2}\int^{R_{1,2}}_{c_0}\frac{dt}{P(t)}=(-1)^{m+1}a_1^{m-1}D_1^{-m+1/2}\frac{\Gamma(1-m)^2}{\Gamma(2(1-m))}\frac{1}{r_{1,2}+1}.
\end{equation}
However $\beta_{r_{1,2}}$ are solutions of (169) and we finaly have
\begin{equation}
\int^{-\rho_1-\frac{\sqrt{D_1}}{a_1}\beta_{r_2}}_{-\rho_1-\frac{\sqrt{D_1}}{a_1}\beta_{r_1}}\frac{f_1(t)}{(a_1t^2+b_1t+c_1)^m}dt=R_2-R_1.
\end{equation}
\\
\textbf{Proof.}\\
Given any $f_1(A)$ and $A_1,A_2$, we have
\begin{equation}
\int^{U_i\left(A_2\right)}_{U_i\left(A_1\right)}\frac{f_1(t)}{(a_1t^2+b_1t+c_1)^m}dt=\pi^2\left(\frac{1}{L'(A_2)^2}-\frac{1}{L'(A_1)^2}\right),
\end{equation}
where $f_1$ and $L$ are related as
\begin{equation}
f_1(A)=\pi^2\left(\frac{d}{dt}\frac{1}{L'(t)^2}\right)_{t=U(A)}.
\end{equation}
But from equation (148) we have
\begin{equation}
\int^{U_i\left(c-\frac{1}{2}\int^{A_2}_{c_0}w'(q)qdA\right)}_{U_i\left(c-\frac{1}{2}\int^{A_1}_{c_0}w'(q)qdA\right)}\frac{f_1(t)}{(a_1t^2+b_1t+c_1)^m}dt=A_2-A_1,
\end{equation}
Assume that $R_{1,2}$ are solutions of
\begin{equation}
c-\frac{1}{2}\int^{R_{1,2}}_{c_0} w'(q)qdA=(-1)^{m+1}a_1^{m-1}D_1^{-m+1/2}\frac{\Gamma(1-m)^2}{\Gamma(2(1-m))}\frac{1}{r_{1,2}+1}
\end{equation}
and $\beta_{r_1},\beta_{r_2}$ are solutions of (169). Then
$$
U_i\left(c-\frac{1}{2}\int^{R_{1,2}}_{c_0} w'(q)qdA\right)=
$$
$$
=U_i\left((-1)^{m+1}a_1^{m-1}D_1^{-m+1/2}\frac{\Gamma(1-m)^2}{\Gamma(2(1-m))}\frac{1}{r_{1,2}+1}\right)=
$$
$$
=-\rho_1-\frac{\sqrt{D_1}}{a_1}\beta_{r_{1,2}}.
$$
Hence we get the proof of the theorem.\\
\\
\textbf{Remarks.} We have 
$$
G(y(A))=1/h_i'(A).
$$ 
Hence
$$
5\int^{y(A)}_{0}\frac{dt}{t\sqrt[6]{t^{-5}-11-t^5}}=h_i(A)\Leftrightarrow F_1(h_i(A))=y(A)\Leftrightarrow
$$
$$
F_1\left(c-\frac{1}{2}\int w'(q)qdA\right)=y(A)\Rightarrow 
$$
$$
G\left(F_1\left(c-\frac{1}{2}\int w'(q)qdA\right)\right)=1/h_i'(A)=-\frac{2}{w'(q)q}\Rightarrow
$$
$$
G(y(A))+2P(A)=0.
$$
Also
$$
G(F_1(A))=-\frac{2}{L(A)}+P^{*}_{0}(A)
$$
$$
f_1(U_i(A))=h'(A)=\pi^2\frac{d}{dA}\frac{1}{L'(A)^2}\Leftrightarrow
$$
\begin{equation}
h(A)=\frac{\pi^2}{L'(A)^2}+l_1,
\end{equation}
where $l_1$ is constant. But differentiating (158) we have
\begin{equation}
L'\left(c-\frac{1}{2}\int w'(q) q dA\right)=\frac{\pi}{\sqrt{A}}.
\end{equation}
Setting in (185) $A\rightarrow c-\frac{1}{2}\int w'(q)qdA$ and using (186),(160), we get
\begin{equation}
h\left(L_i\left(w(q)\right)\right)=A+l_1.
\end{equation}
But $\pi L_i'(w(q))=\sqrt{A}$. Hence
$$
L_i'(w(q))w'(q)q\frac{-\pi}{2\sqrt{A}}=w'(q)q\frac{-\pi}{2\sqrt{A}}\frac{\sqrt{A}}{\pi}=-\frac{w'(q)q}{2}.
$$
Hence from Theorem 14 we get
\begin{equation}
L_i(w(q))=h_i(A)+l_2,
\end{equation}
where $l_2$ is constant. Hence 
$$
h(h_i(A)+l_2)=A+l_1\Rightarrow h_i(A+l_1)=h_i(A)+l_2.\eqno{(188.1)}
$$
Hence the functions $h(A)$, $L_i\left(w\left(e^{-\pi\sqrt{A}}\right)\right)$ are one to one and hence strictly increasing or decreasing. Also their derivatives are periodic. Another intersting thing is that (using (188)):
\begin{equation}
L(A)=w\left(e^{-\pi\sqrt{h(A)-l_1}}\right).
\end{equation}
Also from (185) we have
\begin{equation}
L(A)=\pm \pi \int^{A}_{c^{*}}\frac{dt}{\sqrt{h(t)-l_1}}.
\end{equation}
However if we know $h(A)$ we know by simple inversions the functions $y(A),G(A)$,
$P(A),\int\frac{dA}{P(A)}$ (from relations (133),(156),(148.1),(150)). Hence if we know the expansion
\begin{equation}
G\left(F_1(A)\right)=h'(A)=-\frac{2}{L(A)}+P_0^{*}(A),
\end{equation}
then we can find $f(A)$ from Theorem 16 and then solve (127) with respect to $w$. Hence we have the next\\
\\
\textbf{Theorem 20.}\\
Assume given a function $h(A)$ we can write it in the form 
\begin{equation}
h'(A)=-\frac{2}{L(A)}+P_0^{*}(A),
\end{equation}
where $L(A)$ is solution of the equation
\begin{equation}
-2\pi^2\frac{L''(A)}{L'(A)^3}+\frac{2}{L(A)}=P_0^{*}(A)
\end{equation}
 and $P_0^{*}(A)$ analyitc. Then $f(A)$ is given from (160) and $w(A)$ from $\frac{w(A)}{f(w(A))}=A$ and holds
\begin{equation}
\pm\pi \int^{A}_{c^{*}}\frac{dt}{\sqrt{h(t)-l_1}}=w\left(e^{-\pi\sqrt{h(A)-l_1}}\right)=L(A).
\end{equation}
\\

But
$$
\frac{1}{L'(A)}=-\frac{1}{\pi^2}\log\left(\frac{L(A)}{f(L(A))}\right)\Rightarrow 
$$
$$
-\frac{\pi^2}{L'(A)}=-\pi \sqrt{h(A)-l_1}=\log\left(\frac{L(A)}{f(L(A))}\right)\Rightarrow 
$$
\begin{equation}
e^{-\pi\sqrt{h(A)-l_1}}=\frac{L(A)}{f(L(A))}.
\end{equation}
Set now $Q(A)$ such $L(A)=Q\left(e^{-\pi\sqrt{h(A)-l_1}}\right)$ and $Q(Q(A))=\phi(A)$. Then (here finction $\phi$ must not confused with $\phi$ of [2] and $\phi$ of Section 3 below), we have:
\begin{equation}
A=\frac{Q(A)}{f(Q(A))}\Rightarrow Q\left(A\right)=\frac{\phi(A)}{f\left(\phi(A)\right)}\Rightarrow
L(A)=\frac{\phi\left(e^{-\pi\sqrt{h(A)-l_1}}\right)}{f\left(\phi\left(e^{-\pi \sqrt{h(A)-l_1}}\right)\right)}.
\end{equation}
If we assume that
\begin{equation}
Q(A)=\frac{\phi(A)}{f(\phi(A))}=\phi(\lambda(A)),
\end{equation}
we must have equivalently
$$
Q(A)=\frac{\phi(A)}{f(\phi(A))}\Leftrightarrow Q(Q(A))=\phi(A)\Leftrightarrow
\frac{\phi\left(\frac{\phi(A)}{f(\phi(A))}\right)}{f\left(\phi\left(\frac{\phi(A)}{f(\phi(A))}\right)\right)}=\phi(A)\Leftrightarrow
$$
$$
\frac{\phi(\phi(\lambda(A)))}{f\left(\phi(\phi(\lambda(A)))\right)}=\phi(A).\eqno{(a)}
$$
But from (196) we have
$$
\frac{\phi(\phi(\lambda(A)))}{f(\phi(\phi(\lambda(A))))}=\phi(\lambda(\phi(\lambda(A)))).\eqno{(b)}
$$
Hence from $(a),(b)$, we must have
$$
\phi(\lambda(\phi(\lambda(A))))=\phi(A)\Leftrightarrow \lambda(\phi(\lambda(A)))=A.
$$
Hence
$$
L(A)=\lambda^{(-1)}\left(e^{-\pi\sqrt{h(A)-l_1}}\right).
$$
Hence from (189) we must have 
$$
\lambda^{(-1)}(A)=\frac{\phi(A)}{f(\phi(A))}=w_i(\phi(A))\Leftrightarrow
$$
$$
w(A)=\phi(\lambda(A))=\lambda^{(-1)}(A)=\frac{\phi(A)}{f(\phi(A))}\Leftrightarrow
$$
\begin{equation}
w(w(A))=\phi(A)\textrm{, }\lambda(A)=\frac{A}{f(A)}.
\end{equation}
Also
$$
\lambda\left(\lambda(A)\right)=\phi^{(-1)}(A)\Leftrightarrow \frac{A}{f(A)f\left(\frac{A}{f(A)}\right)}=\phi^{(-1)}(A)
$$
\\
\textbf{Theorem 21.}\\
We have
\begin{equation}
y(A)=F_1\left(-\frac{1}{2}\int^{A}_{c}\frac{dt}{P(t)}\right),
\end{equation}
\begin{equation}
G(y(A))+2P(A)=0
\end{equation}
and
\begin{equation}
\int^{-1/2\int^{A}_{c_1}1/P(t)dt}_{c}G\left(F_1(t)\right)dt=A.
\end{equation}
\\
\textbf{Proof.}\\
We have $-2 P(A) h_i'(A)=1$ and $h_i'(A)=1/G(y(A))$. Hence 
$$
G(y(A))=-2 P(A).
$$
Also $F_1(h_i(A))=y(A)$. Hence
$$
y(A)=F_1\left(-\frac{1}{2}\int^{A}_{c}\frac{dt}{P(t)}\right).
$$

\section{Solving polynomial equations}

An interesting case of functions are the Lambert functions defined as
\begin{equation}
\phi(x)=q^2\sum^{\infty}_{n=1}\frac{A_nq^n}{1-q^n}\textrm{, }q=e^{-\pi\sqrt{x}}
\end{equation} 
where $A_n=\sum_{d|n}a_d\mu(n/d)$ and $a_n$ is arithmetic $T$-periodic function. Then we can write 
$$
\phi(x)=q^2\sum^{\infty}_{n=1}\frac{A_nq^n}{1-q^n}=q^2\sum^{\infty}_{n=1}\left(\sum_{d|n}A_d\right)q^n=q^2(1-q^T)^{-1}\sum^{T}_{n=1}\left(\sum_{d|n}A_d\right)q^n=
$$
\begin{equation}
=q^2(1-q^T)^{-1}\sum^{T}_{n=1}a_nq^n.
\end{equation}
For example if $a_n=\sqrt{2}\cos(\pi n/4)$, then $T=8$ and 
$$
\sqrt{2}\sum^{\infty}_{n=1}\sum_{d|n}\cos(\pi d/4)\mu(n/d)\frac{q^n}{1-q^n}=\frac{q-q^3-\sqrt{2}q^4-q^5+q^7+\sqrt{2}q^8}{1-q^8}
$$
Hence the series 
\begin{equation}
\phi(x)=q^2\sum^{\infty}_{n=1}\frac{A_nq^n}{1-q^n}\textrm{, }q=e^{-\pi\sqrt{x}}\textrm{, }x>0
\end{equation}
is a rational function of $q=e^{-\pi\sqrt{x}}$. Also holds
\begin{equation}
\phi(x)=q^2\sum^{\infty}_{n=1}q^na_n,
\end{equation}
where $a_n=\sum_{d|n}A_d$.\\
 
An analysis of how we find functions like $\phi(x)$ is given in [2]. For example let $G(x)$ be such that  
\begin{equation}
\phi^{(-1)}(k_x)=m^{(-1)}_G(x)=\pi\int^{+\infty}_{\sqrt{x}}\eta(it/2)^4G(R(e^{-\pi t}))dt,
\end{equation}
where 
\begin{equation}
\eta(z)=q^{1/24}\prod^{\infty}_{n=1}(1-q^n)\textrm{, }q=e^{2\pi i z}\textrm{, }Im(z)>0
\end{equation}
and 
\begin{equation}
R(x)=\frac{x^{1/5}}{1+}\frac{x^1}{1+}\frac{x^2}{1+}\frac{x^3}{1+}\ldots\textrm{, }|x|<1,
\end{equation}
are the Dedekind eta function and the Rogers-Ramanujan continued fraction respectively.
Consequently the function $y$ can be found as
\begin{equation}
y(x)=R\left(e^{-\pi\sqrt{k_i(\phi(x))}}\right).
\end{equation}
Hence with the notation of [2] relation (59), (here Lambert's function $\phi(x)$ is not to be confused with the notation of $\phi(x)$ of [2] relation (14)), we have:  
\begin{equation}
\phi(x)=s(x)
\end{equation}
and thus
\begin{equation}
G(x)=\frac{G_0(x)}{\sigma(F_i(x))},
\end{equation}
where 
\begin{equation}
\sigma(x)=\frac{1}{s_i'(x)}=\frac{1}{\phi^{(-1)}{ }{'}(x)}.
\end{equation}
Hence
\begin{equation}
G(x)=G_0(x)\phi^{(-1)}{ }{'}(F_i(x)).
\end{equation}
Now we have in general
\begin{equation}
F_1\left(h_i'(x)\right)=y(x).
\end{equation}
Hence
$$
G\left(F_1\left(h_i'(x)\right)\right)=G(y(x))
$$
and for $X(x)$ it holds that 
\begin{equation}
X(A)=s_i(A)=\phi^{(-1)}(x)
\end{equation}
and
\begin{equation}
G\left(F_1\left(-\frac{1}{2}\int^{x}_{c}\frac{dt}{P(t)}\right)\right)=-2P(x)
\end{equation}
and
\begin{equation}
G(y(x))=-2P(x)\textrm{ and }-2h_i'(x)P(x)=1.
\end{equation}
Hence from [3] Theorem 12 we get the next\\ 
\\
\textbf{Theorem 22.}
\begin{equation}
P(x)=-\frac{\left(\phi(x)\sqrt{1-\phi(x)^2}\right)^{2/3}}{2\sqrt[3]{2}\phi'(x)}
\end{equation}
and
\begin{equation}
h_i(x)=c-\frac{1}{2}\int^{x}_{c_0}\frac{dt}{P(t)}=\frac{1}{\sqrt[3]{4}}B_0\left(\phi(x)^2;\frac{1}{6},\frac{2}{3}\right),
\end{equation}
where $B_0(x;a,b)=\int^{x}_{0}t^{a-1}(1-t)^{b-1}dt$ is the incomplete Beta function. Also
$$
F(x)=-2h_i(x)+c'\textrm{ and }y(x)=R\left(e^{-\pi \sqrt{k_i(\phi(x))}}\right).\eqno{(219.1)}
$$
\\

Someone can also easily see that
\begin{equation}
\int\frac{dx}{P(x)}=-\frac{2}{\sqrt[3]{4}}B_0\left(\phi(x)^2;\frac{1}{6},\frac{2}{3}\right)+c,
\end{equation}
\begin{equation}
P\left(X(x)\right)=-\frac{\left(x\sqrt{1-x^2}\right)^{2/3}}{2\sqrt[3]{2}}\phi^{(-1)}{'}(x)
\end{equation}
and
\begin{equation}
Y(x)=\phi^{(-1)}(k_x).
\end{equation}
There exist $\alpha>0$ such that for all $x\in[0,\alpha)$, we have ($\phi(0)=0$),
\begin{equation}
F(x)=\frac{1}{\sqrt[3]{4}}B_0\left(\phi(x)^2;\frac{1}{6},\frac{2}{3}\right)=h_i(x)
\end{equation}
and hence
$$
F\left(\phi^{(-1)}\left(x^2\right)\right)+F\left(\phi^{(-1)}\left[\left(\frac{1-x}{1+x}\right)^2\right]\right)=\frac{\sqrt{3}\Gamma\left(\frac{1}{3}\right)^3}{2\pi\sqrt[3]{2}}\Leftrightarrow
$$
$$
F\left(\phi^{(-1)}\left(x\right)\right)+F\left(\phi^{(-1)}\left[\left(\frac{1-\sqrt{x}}{1+\sqrt{x}}\right)^2\right]\right)=\frac{\sqrt{3}\Gamma\left(\frac{1}{3}\right)^3}{2\pi\sqrt[3]{2}}\Leftrightarrow
$$
\begin{equation}
F\left(x\right)+F\left(\phi^{(-1)}\left[\left(\frac{1-\sqrt{\phi(x)}}{1+\sqrt{\phi(x)}}\right)^2\right]\right)=\frac{\sqrt{3}\Gamma\left(\frac{1}{3}\right)^3}{2\pi\sqrt[3]{2}}.
\end{equation}
Also $-\frac{1}{2}F'(x)=-\frac{1}{2P(x)}=h_i'(x)$. \\
\\
\textbf{Theorem 22.1} (The complex analog.)\\
If we define $m^{*}(z)$ as
\begin{equation}
m^{*}(z)=\left(\frac{\theta_2\left(e^{i\pi z}\right)}{\theta_3\left(e^{i\pi z}\right)}\right)^2\textrm{, }Im(z)>0,
\end{equation}
then 
\begin{equation}
m^{*}\left(-\frac{1}{z}\right)=\sqrt{1-m^{*}(z)^2}\textrm{, }m^{*}(z+1)=\frac{m^{*}(z)}{\sqrt{m^{*}(z)^2-1}}
\end{equation}
and
\begin{equation}
m^{*}(z+2)=m^{*}(z).
\end{equation}
Also
\begin{equation}
\phi(A)=s(A),
\end{equation}
\begin{equation}
\phi^{(-1)}\left(m^{*}(2z)\right)=s^{(-1)}\left(m^{*}(2z)\right)=Y(z)
\end{equation}
and $Y(z)$ is Hauptmodul satisfying
\begin{equation}
Y'(z)+4\pi i\cdot \eta(z)^4P(Y(z))=0,
\end{equation}
\begin{equation}
F\left(Y(z)\right)+F\left(Y\left(-\frac{1}{z}\right)\right)=\frac{\sqrt{3}\Gamma\left(\frac{1}{3}\right)^3}{2\pi\sqrt[3]{2}},
\end{equation}
where
\begin{equation}
F(z)=\frac{1}{\sqrt[3]{4}}B_0\left(\phi(z)^2,\frac{1}{6},\frac{2}{3}\right)=h_i(z).
\end{equation}
Also
\begin{equation}
5\int^{y(A)}_{0}\frac{G(t)}{t\sqrt[6]{t^{-5}-11-t^5}}dt=A,
\end{equation}
where
\begin{equation}
G(A)=\frac{G_0(A)}{\sigma\left(F_0^{(-1)}(A)\right)}=2^{-1/3}\frac{\left(F^{(-1)}_0(A)\sqrt{1-F^{(-1)}_0(A)^2}\right)^{2/3}}{\sigma\left(F^{(-1)}_0(A)\right)}
\end{equation}
and $F_0(A)=R\left(e^{i\pi {m^{*}}^{(-1)}(A)}\right)$, where ${m^{*}}^{(-1)}(A)$ is the inverse of $m^{*}(A)$, in the sense
\begin{equation}
{m^{*}}^{(-1)}(A)=i\cdot\frac{{}_2F_1\left(\frac{1}{2},\frac{1}{2};1;1-A^2\right)}{{}_2F_1\left(\frac{1}{2},\frac{1}{2};1;A^2\right)},
\end{equation}
where $-\frac{1}{2}<Re(z)\leq\frac{1}{2}$, $Im(z)>0$.
\begin{equation}
\phi(z)=\sum^{\infty}_{n=1}\frac{A_nq^n}{1-q^n}\textrm{, }q=e^{i\pi z}\textrm{, }Im(z)>0,
\end{equation}
\begin{equation}
\sigma(A)=\frac{1}{\phi^{(-1)}{'}(A)},
\end{equation}
\begin{equation}
m_{G}^{(-1)}(A):=2\phi^{(-1)}(m^{*}(A))=2i \pi\int^{A}_{+i\infty}\eta(t/2)^4G\left(R\left(e^{i\pi t}\right)\right)dt\textrm{, }Im(A)>0.
\end{equation}
The solution of (233) is
\begin{equation}
y(A)=R\left(e^{i\pi {m^{*}}^{(-1)}(\phi(A))}\right)\Leftrightarrow y\left(Y(A)\right)=R\left(e(A)\right).
\end{equation}
\\
\textbf{Theorem 22.2} (In the complex analog)\\
If we assume the problems (233), (230), then the functions $G(A)$, $P(A)$ are related with equation
$$
G\left(F_1\left(\frac{c}{2\pi i}-\int^{P_i(A)}_{i\infty}\frac{dt}{P(t)}\right)\right)=-A.\eqno{(239.1)}
$$
Also if
$$
F(A)=\int^{A}_{Y(i\infty)}\frac{dt}{P(t)},
$$
then
$$
m_G^{(-1)}\left(-\frac{2}{A}\right)=F^{(-1)}\left(\frac{\sqrt{3}\Gamma\left(\frac{1}{3}\right)^3}{2\pi\sqrt[3]{2}}-F\left(m_G^{(-1)}(2A)\right)\right)\eqno{(239.2)}
$$
and
$$
Y(A)=m_G^{(-1)}\left(2A\right)=2\pi i\int^{2A}_{+i\infty}\eta(t/2)^4G\left(R(e^{\pi i t})\right)dt,\eqno{(239.3)}
$$
$$
G\left(R\left(q_A\right)\right)=-P(Y(A))\textrm{, }q_A=e(A)\textrm{, }Im(A)>0.\eqno{(239.4)}
$$
Moreover if 
$$
\frac{1}{P(z)}=\sum^{\infty}_{n=1}a_nq^n\textrm{, }q=e(z)\textrm{, }Im(z)>0,
$$
then $G(y(z))=-P(z)$, $y(z)=R\left(e^{i\pi m_{G}(z)}\right)$ and exists constant $C$ such that:
$$
\prod^{\infty}_{n=1}\left(1-q^n\right)^{-1/n\sum_{d|n}a_d\mu(n/d)}=C\exp\left(8\pi^2\int^{m_G(z)/2}_{+i\infty}\eta(t)^4dt\right)\eqno{(239.5)}
$$
and
$$
w(q_A)=8\pi^2 \int^{m_{G}(A)/2}_{+i\infty}\eta(t)^4dt.\eqno{(239.6)}
$$
Hence
$$
e^{w(q)}=C'\prod^{\infty}_{n=1}\left(1-q^n\right)^{-1/n\sum_{d|n}a_d\mu(n/d)}.\eqno{(239.7)}
$$
\textbf{Remarks.}\\
From equation (239.1) we can write
$$
P^{(-1)}{'}(A)=-A\frac{d}{dA}\left[F_1^{(-1)}\left(G^{(-1)}\left(-A\right)\right)\right].
$$
\\
\textbf{Proof.}
$$
y(A)=R\left(e^{i\pi m^{{*}{(-1)}}(\phi(A))}\right)=R\left(e^{i\pi m_G(A)}\right)\Rightarrow 
$$
$$
y\left(Y(A)\right)=R\left(e^{i\pi m_{G}(Y(A))}\right)=R\left(e^{2\pi i A}\right)\Rightarrow m_G\left(Y(A)\right)=2A\Rightarrow
$$
$$
Y(A)=m_G^{(-1)}\left(2A\right)=2\pi i\int^{2A}_{+i\infty}\eta(t/2)^4G\left(R(e^{\pi i t})\right)dt.
$$
Hence
$$
F\left(m_G^{(-1)}\left(2A\right)\right)+F\left(m_G^{(-1)}\left(-\frac{2}{A}\right)\right)=\frac{\sqrt{3}\Gamma\left(\frac{1}{3}\right)^3}{2\pi\sqrt[3]{2}}.
$$
\[
\]

Now we define $a^{*}_n$, $b^{*}_n$ such that
\begin{equation}
\frac{\sqrt{x}}{P(x)}=\sum^{\infty}_{n=1}a^{*}_ne^{-\pi n\sqrt{x}}=-\frac{2\sqrt[3]{2}\phi'(x)\sqrt{x}}{\left(\phi(x)\sqrt{1-\phi(x)^2}\right)^{2/3}}
\end{equation}
and
\begin{equation}
\frac{\sum^{\infty}_{n=1}a^{*}_ne^{-\pi n\sqrt{x}}}{\sum^{\infty}_{n=1}b^{*}_ne^{-\pi n \sqrt{x}}}=\sqrt{x}.
\end{equation} 
Then we set $c^{(1)}_n:=\frac{a^{*}_n}{n}$ and $c^{(2)}_n:=\frac{b^{*}_n}{n}$ and assume that
\begin{equation}
w_1(q)=\sum^{\infty}_{n=1}c^{(1)}_nq^n\textrm{, }w_2(q)=\sum^{\infty}_{n=1}c_n^{(2)}q^n,
\end{equation}
are the solutions of the equations
\begin{equation}
\frac{w_1(q)}{f_1(w_1(q))}=q\textrm{, }\frac{w_2(q)}{f_2(w_2(q))}=q\textrm{, resp.},
\end{equation}
where $f_1(x)$,$f_2(x)$ are functions such that $f_1(0)$,$f_2(0)\neq 0$ and analytic arround the origin, we get integrating (240):
$$
\int\left(\sum^{\infty}_{n=1}a_n^{*}\frac{e^{-\pi n \sqrt{x}}}{\sqrt{x}}\right)dx=-\frac{2}{\pi}\frac{\pi}{\sqrt[3]{4}} B_0\left(\phi(x)^2;\frac{1}{6},\frac{2}{3}\right)+c\Leftrightarrow
$$
$$
\frac{-2}{\pi}w_1(q)=-\frac{2}{\pi}\sum^{\infty}_{n=1}c^{(1)}_nq^n=-\frac{2}{\pi}\frac{\pi}{\sqrt[3]{4}}B_0\left(\phi(x)^2;\frac{1}{6},\frac{2}{3}\right)\Leftrightarrow
$$
$$
w_1(q)=\sum^{\infty}_{n=1}c^{(1)}_nq^n=\frac{\pi}{\sqrt[3]{4}}B_0\left(\phi(x)^2;\frac{1}{6},\frac{2}{3}\right).
$$
Note that
$$
\phi(0)=(+\infty) \left(\sum^{T}_{k=1}a_k\right)\textrm{, }\phi(+\infty)=0.
$$
Hence $\phi^{(-1)}(0)=+\infty\Rightarrow X(0)=+\infty$ and
$$
F(x)=\frac{1}{\sqrt[3]{4}}B_0\left(\phi(x)^2;\frac{1}{6},\frac{2}{3}\right).
$$
Hence when $q=e^{-\pi\sqrt{x}}$, $x>0$, we have
\begin{equation}
w_1(q)=\frac{\pi}{\sqrt[3]{4}}B_0\left(\phi(x)^2;\frac{1}{6},\frac{2}{3}\right)=\pi F(x)=\pi h_i(x).
\end{equation}
Also
$$
\int\left(\sum^{\infty}_{n=1}b_n^{*}\frac{e^{-\pi n \sqrt{x}}}{\sqrt{x}}\right)dx=-2\sqrt[3]{2}\int^{x}_{c}\frac{\phi'(t)}{\sqrt{t}\left(\phi(t)\sqrt{1-\phi(t)^2}\right)^{2/3}}dt\Leftrightarrow
$$
$$
\frac{-2}{\pi}w_2(q)=-\frac{2}{\pi}\sum^{\infty}_{n=1}c^{(2)}_nq^n=-2\sqrt[3]{2}\int^{x}_{c}\frac{\phi'(t)}{\sqrt{t}\left(\phi(t)\sqrt{1-\phi(t)^2}\right)^{2/3}}dt\Leftrightarrow
$$
\begin{equation}
w_2(q)=\pi\sqrt[3]{2}\int^{x}_{+\infty}\frac{\phi'(t)}{\sqrt{t}\left(\phi(t)\sqrt{1-\phi(t)^2}\right)^{2/3}}dt.
\end{equation}
However from (131) it holds
$$
\int\frac{1}{qP(A)}dq=w(q)+c\Leftrightarrow w(q)=\int\frac{1}{qP(x)}q\frac{-\pi}{2\sqrt{x}}dx+c=
$$
$$
=-\frac{\pi}{2}\int\frac{1}{P(x)\sqrt{x}}dx=-\frac{\pi}{2}(-2\sqrt[3]{2})\int\frac{\phi'(x)}{\sqrt{x}\left(\phi(x)\sqrt{1-\phi(x)^2}\right)^{2/3}}dx\Leftrightarrow
$$
$$
w(q)=w_2(q)+c.
$$
Hence we get the next\\
\\
\textbf{Theorem 23.}\\
If $q=e^{-\pi\sqrt{x}}$, $x>0$, then
\begin{equation}
w_1(q)=\sum^{\infty}_{n=1}\frac{a_n^{*}}{n}e^{-\pi n\sqrt{x}}=\frac{\pi}{\sqrt[3]{4}}B_0\left(\phi(x)^2;\frac{1}{6},\frac{2}{3}\right)=\pi F(x)=\pi h_i(x)
\end{equation}
and
\begin{equation}
w(q)=w_2(q)=\sum^{\infty}_{n=1}\frac{b_n^{*}}{n}e^{-\pi n \sqrt{x}}=\pi\sqrt[3]{2}\int^{x}_{+\infty}\frac{\phi'(t)}{\sqrt{t}\left(\phi(t)\sqrt{1-\phi(t)^2}\right)^{2/3}}dt,
\end{equation}
where
\begin{equation}
\frac{\sum^{\infty}_{n=1}a_n^{*}e^{-\pi n \sqrt{x}}}{\sum^{\infty}_{n=1}b_n^{*}e^{-\pi n \sqrt{x}}}=\sqrt{x}.
\end{equation}
\textbf{Remarks.} Knowing $a_n^{*}$ from (246), we can find $b_n^{*}$ from (248) and hence $w(q)$ from (247). By this way we can evaluate $F(x),h_i(x),L(x),y(x),P(x),Y(x)$ (for $L(x)$ see relation (263) below), and $f(x)$ with simple inversion. From (200) we get 
$$
G\left(R\left(e^{-\pi\sqrt{k_i(\phi(x))}}\right)\right)=\frac{\left(\phi(x)\sqrt{1-\phi(x)^2}\right)^{2/3}}{\sqrt[3]{2}\phi'(x)}\Leftrightarrow
$$
$$
G\left(F_0(x)\right)=\frac{\left(x\sqrt{1-x^2}\right)^{2/3}}{\sqrt[3]{2}}\phi_i'(x)
$$
\\
\textbf{Lemma 1.}\\
If $\frac{w_1(q)}{f_1\left(w_1(q)\right)}=q$, with $f_1,w_1$ as above, then
\begin{equation}
\frac{\pi}{\sqrt[3]{4}}B_0\left(\phi\left(\pi^{-2}\log^2(x/f_1(x))\right)^2;\frac{1}{6},\frac{2}{3}\right)=x.
\end{equation}
\\
\textbf{Proof.}\\
$$
w_1\left(e^{-\pi \sqrt{x}}\right)=\frac{\pi}{\sqrt[3]{4}}B_0\left(\phi(x)^2;\frac{1}{6},\frac{2}{3}\right)\Leftrightarrow
$$
$$
w_1\left(\frac{1}{x}\right)=\frac{\pi}{\sqrt[3]{4}}B_0\left(\phi\left(\pi^{-2}\log^2 x\right)^2;\frac{1}{6},\frac{2}{3}\right).
$$
Hence $w_1(1/x)=w_1(x)$. Also
$$
x=\frac{\pi}{\sqrt[3]{4}}B_0\left(\phi\left(\pi^{-2}\log^2\left(x/f_1(x)\right)\right)^2;\frac{1}{6},\frac{2}{3}\right).
$$
\\

In the same way as above\\
\\
\textbf{Lemma 2.}\\
If $\frac{w(q)}{f\left(w(q)\right)}=q$, with $f(x)$ analytic at the origin and $f(0)\neq 0$, then
\begin{equation}
\pi\sqrt[3]{2}\int^{\pi^{-2}\log^2\left(x/f(x)\right)}_{+\infty}\frac{\phi'(t)}{\sqrt{t}\left(\phi(t)\sqrt{1-\phi(t)^2}\right)^{2/3}}dt=x.
\end{equation}
\\

Also if we define the function $m(x)$ such that
\begin{equation}
\pi\int^{+\infty}_{\sqrt{m(x)}}\eta\left(it/2\right)^4dt=x,
\end{equation}
then
\begin{equation}
\frac{1}{\sqrt[3]{4}}B_0\left(k(m(x))^2;\frac{1}{6},\frac{2}{3}\right)=x
\end{equation}
and
$$
F_1(A)=R\left(e^{-\pi\sqrt{m(A)}}\right).\eqno{(252.1)}
$$
Hence we have the next\\
\\
\textbf{Lemma 3.}\\
If $f_1(x)$ is analytic arround 0 and $f_1(0)\neq 0$, then $\phi_i(0)=+\infty$
\begin{equation}
f_1(x)=x\exp\left[\pi\sqrt{\phi_i\left(k\left(m\left(\frac{x}{\pi}\right)\right)\right)}\right].
\end{equation}
\\
\textbf{Proof.}\\
From
$$
w_1\left(e^{-\pi\sqrt{x}}\right)=\frac{\pi}{\sqrt[3]{4}}B_0\left(\phi(x)^2;\frac{1}{6},\frac{2}{3}\right)\Leftrightarrow
w_1\left(e^{-\pi\sqrt{\phi_i(x)}}\right)=\frac{\pi}{\sqrt[3]{4}}B_0\left(x^2;\frac{1}{6},\frac{2}{3}\right)\Leftrightarrow
$$
$$
w_1\left(e^{-\pi\sqrt{\phi_i(k(x))}}\right)=\frac{\pi}{\sqrt[3]{4}}B_0\left(k(x)^2;\frac{1}{6},\frac{2}{3}\right)\Leftrightarrow w_1\left(e^{-\pi\sqrt{\phi_i\left(k\left(m\left(\frac{x}{\pi}\right)\right)\right)}}\right)=x\Leftrightarrow
$$
$$
f_1(x)=xe^{\pi\sqrt{\phi_i\left(k\left(m\left(x/\pi\right)\right)\right)}},
$$
we get the result.\\
\\

Now from Theorem 16 we have
\begin{equation}
w\left(e^{-\pi^2 L_i'(x)}\right)=x.
\end{equation}
Also from (188) we have
\begin{equation}
L\left(\frac{1}{\sqrt[3]{4}}B_0\left(\phi(x)^2;\frac{1}{6},\frac{2}{3}\right)+l_2\right)=\pi\sqrt[3]{2}\int^{x}_{+\infty}\frac{\phi'(t)}{\sqrt{t}\left(\phi(t)\sqrt{1-\phi(t)^2}\right)^{2/3}}dt
\end{equation}
and (from (244),(240))
\begin{equation}
L^{(-1)}{ }{'}{ }^{(-1)}(x)=\pi \sqrt[3]{2}\int^{\pi^2 x^2}_{+\infty}\frac{\phi'(t)}{\sqrt{t}\left(\phi(t)\sqrt{1-\phi(t)^2}\right)^{2/3}}dt.
\end{equation}
Hence differentiating the relation (256) we get
$$
2\pi^2\sqrt[3]{2}\frac{\phi'(x^2\pi^2)}{\left(\phi(x^2\pi^2)\sqrt{1-\phi\left(x^2\pi^2\right)^2}\right)^{2/3}}=L^{(-1)}{ }{'}{ }^{(-1)}{ }{'}\left(x\right).
$$
Hence\\
\\
\textbf{Lemma 4.}
\begin{equation}
L^{(-1)}{ }{'}{ }^{(-1)}{ }{'}\left(x\right)=-\frac{\pi^2}{P\left(x^2\pi^2\right)}.
\end{equation}
\\

Continuing from (255) we get
\begin{equation}
L\left(\frac{1}{\sqrt[3]{4}}B_0\left(\phi\left(\pi^2 x^2\right)^2;\frac{1}{6},\frac{2}{3}\right)\right)=L^{(-1)}{ }{'}{ }^{(-1)}(x).
\end{equation}
Hence
\begin{equation}
\frac{1}{L'\left(\frac{1}{\sqrt[3]{4}}B_0\left(\phi(x)^2;\frac{1}{6},\frac{2}{3}\right)\right)}=\frac{\sqrt{x}}{\pi}.
\end{equation}
From this we have the next\\
\\
\textbf{Theorem 24.}
\begin{equation}
L'\left(\frac{1}{\sqrt[3]{4}}B_0\left(\phi(x)^2;\frac{1}{6},\frac{2}{3}\right)\right)=\frac{\pi}{\sqrt{x}}
\end{equation}
and
\begin{equation}
\pi\sqrt[3]{2}\int^{x}_{+\infty}\frac{\phi'(t)}{\sqrt{t}\left(\phi(t)\sqrt{1-\phi(t)^2}\right)^{2/3}}dt=L\left(\frac{1}{\sqrt[3]{4}}B_0\left(\phi(x)^2;\frac{1}{6},\frac{2}{3}\right)+l_2\right).
\end{equation}
\\

Relation (254) is equivalent to
\begin{equation}
L\left(\pi^{-1}w_1(q)+l_2\right)=w(q).
\end{equation}
\\

Continuing we have from (255), (using the function $m(x)$):
$$
L'\left(\frac{1}{\sqrt[3]{4}}B\left(k(x)^2;\frac{1}{6},\frac{2}{3}\right)+l_2\right)=\frac{\pi}{\sqrt{\phi^{(-1)}(k(x))}}\Rightarrow
$$
$$
L'\left(\pi\int^{+\infty}_{\sqrt{x}}\eta(it/2)^4dt+l_2\right)=\frac{\pi}{\sqrt{\phi^{(-1)}(k(x))}}\Rightarrow
$$
$$
L'\left(x+l_2\right)=\frac{\pi}{\sqrt{\phi^{(-1)}(k(m(x)))}}.
$$
Hence\\
\\
\textbf{Theorem 25.}
\begin{equation}
L(x)=\pi\int^{x}_{c}\frac{1}{\sqrt{\phi^{(-1)}(k(m(t)))-l_1}}dt=w\left(e^{-\pi\sqrt{\phi^{(-1)}(k(m(x)))-l_1}}\right)+c'_1
\end{equation}
and
\begin{equation}
h(x)=\phi^{(-1)}\left(k(m(x))\right).
\end{equation}
\\
$$
x=\pi\int^{L_i(x)-l_2}_{c}\frac{1}{\sqrt{\phi^{(-1)}(k(m(t)))}}dt\Rightarrow 1=\pi\frac{L^{(-1)}{ }{'}(x)}{\sqrt{\phi^{(-1)}(k(m(L_i(x))))}}\Rightarrow
$$
$$
L^{(-1)}{ }{'}(x)=\pi^{-1}\sqrt{\phi^{(-1)}\left(k\left(m\left(L^{(-1)}(x)\right)\right)\right)}.
$$

\section{Lagrange inversion in the complex analog}

Using equations (1)-(6) we can easily show that if $q=e^{-\pi\sqrt{x}}$, $x>0$ and $g(q)$ is analytic arround 0 with $g(0)=1$:
$$
g(q)=e^{w(q)}=1+\sum^{\infty}_{n=1}\frac{q^n}{n!}\left(\frac{D}{Dh}\right)^{n-1}\left(e^hf(h)^n\right)_{h=0}=
$$
\begin{equation}
=\prod^{\infty}_{n=1}\left(1-q^n\right)^{-\frac{1}{n}\sum_{d|n}\frac{\mu(n/d)}{\Gamma(d)}\left[\left(\frac{D}{Dh}\right)^{d-1}\left(f(h)^d\right)\right]_{h=0}},
\end{equation}
where $f(x)$ is defined as follows: It holds
\begin{equation}
w(x)=\log\left(g(x)\right).
\end{equation}
If $g_1(x)$ is the inverse of $\log(g(x))$, then
\begin{equation}
f(x)=\frac{x}{g_1(x)},
\end{equation}
where $g_1(x)$ is analytic around 0 and have simple root at $x=0$ i.e. $g_1(0)=0$ and $g_1'(0)\neq 0$.\\
\\
\textbf{Example.}\\
Asume that
$$
g(q)=\sqrt{1+q+q^2},
$$
then set $w(x)=\log\left(\sqrt{1+x+x^2}\right)$. Solving $w(x)=y$, we get
$$
x=w^{(-1)}(y)=\frac{1}{2}\left(-1+\sqrt{-3+4e^{2y}}\right).
$$
Hence 
$$
f(y)=\frac{2y}{-1+\sqrt{-3+4e^{2y}}}.
$$
Set
\begin{equation}
A_0(n)=\frac{1}{n}\sum_{d|n}\frac{\mu(n/d)}{\Gamma(d)}\left[\left(\frac{D}{Dh}\right)^{d-1}\left(f(h)^d\right)\right]_{h=0}.
\end{equation} 
Then
\begin{equation}
g(q)=\prod^{\infty}_{n=1}\left(1-q^n\right)^{-A_0(n)}.
\end{equation}
Here (in our example) we get $A_0(1)=\frac{1}{2}$, $A_0(3)=-\frac{1}{2}$ and $A_0(n)=0$ for $n\neq 1,3$. Hence
\begin{equation}
\sqrt{1+q+q^2}=(1-q)^{-1/2}(1-q^3)^{1/2}.
\end{equation}
\\

Hence we have the next:\\
\\
\textbf{Theorem 25.1}\\
For the analytic function $g(q)$ holds
\begin{equation}
g(q)=\prod^{\infty}_{n=1}\left(1-q^n\right)^{-1/n\sum_{d|n}\frac{\mu(n/d)}{\Gamma(d)}\left[\left(\frac{D}{Dh}\right)^{d-1}\left(h/g^{(-1)}(e^h)\right)^d\right]_{h=0}}.
\end{equation}
\\
\textbf{Corollary 1.}\\
For every function $g$ analytic in the origin, with $g(0)=1$, $g'(0)\neq 0$, we have 
\begin{equation}
\frac{1}{n!}\left[\left(\frac{D}{Dh}\right)^{n-1}\left(e^h\left(\frac{h}{g^{(-1)}\left(e^h\right)}\right)^n\right)\right]_{h=0}=\frac{g^{(n)}(0)}{n!}\textrm{, }n=1,2,\ldots
\end{equation}
\\
\textbf{Corollary 2.}\\
Assume that $f(z),g(z)$ are functions such that
\begin{equation}
f(z)=1+c_1z+c_2z^2+\ldots\textrm{, }c_1\neq 0
\end{equation}
and
\begin{equation}
g^{(-1)}(z)=\log(z) \left(1+c_1\log z+c_2\log^2 z+\ldots\right)^{-1}.
\end{equation}
Then the equation
\begin{equation}
\frac{w(q)}{f(w(q))}=q,
\end{equation}
have solution $w(q)$ such that
\begin{equation}
e^{w(q)}=g(q)\textrm{, }|q|<C,
\end{equation}
where $C$ is suitable constant.\\
\\

Assume the equation
\begin{equation}
\frac{w(q)}{f(w(q))}=q.
\end{equation}
Then its solution is the function 
\begin{equation}
w(q)=\sum^{\infty}_{n=1}c_nq^n\textrm{, }q=e(z)\textrm{, }Im(z)>0,
\end{equation}
where
\begin{equation}
c_n=\frac{1}{n!}\left[\left(\frac{D}{Dh}\right)^{n-1}f(h)^n\right]_{h=0}.
\end{equation}
Also holds
\begin{equation}
\exp\left(w(q)\right)=\prod^{\infty}_{n=1}\left(1-q^n\right)^{-1/n \sum_{d|n}dc_d\mu(n/d)}.
\end{equation}
Hence if
\begin{equation}
\frac{1}{P(z)}=\sum^{\infty}_{n=1}nc_nq^n=\frac{1}{2\pi i}\frac{d}{dz}w(q)=w'(q)q.
\end{equation}
If we invert (277) we find 
\begin{equation}
w^{(-1)}(q)=\frac{q}{f(q)}.
\end{equation}
Hence 
\begin{equation}
f(A)=\frac{A}{w^{(-1)}(A)}.
\end{equation}
From (281) with integration we get ($w(0)=0$):
$$
\int^{z}_{+i\infty}\frac{1}{P(t)}dt=\frac{1}{2\pi i}w(q)\Rightarrow 2\pi i\int^{Y(z)}_{+i\infty}\frac{dt}{P(t)}=w\left(e\left(Y(z)\right)\right)=
$$
$$
=\frac{4\pi}{i}\frac{1}{\sqrt[3]{4}}B\left(m^{*}(2z)^2;\frac{1}{6},\frac{2}{3}\right).
$$
Hence from
$$
\frac{w(e(Y(z)))}{f(w(e(Y(z))))}=e(Y(z)),
$$
it follows that
\begin{equation}
e(Y(z))=\frac{\frac{4\pi}{i\sqrt[3]{4}}B_0\left(m^{*}(2z)^2;\frac{1}{6},\frac{2}{3}\right)}{f\left(\frac{4\pi}{i\sqrt[3]{4}}B_0\left(m^{*}(2z)^2;\frac{1}{6},\frac{2}{3}\right)\right)}.
\end{equation}
Hence assuming the logarithm in the sence if $w=|w|(\cos\phi+i\sin\phi)$, then
\begin{equation}
e^z=w\Leftrightarrow z=\log(w)=\log|w|+i\phi,
\end{equation}
we can write
\begin{equation}
Y(z)=\frac{1}{2\pi i}\log\left(\frac{\frac{4\pi}{i\sqrt[3]{4}}B_0\left(m^{*}(2z)^2;\frac{1}{6},\frac{2}{3}\right)}{f\left(\frac{4\pi}{i\sqrt[3]{4}}B_0\left(m^{*}(2z)^2;\frac{1}{6},\frac{2}{3}\right)\right)}\right).
\end{equation}
But $y(Y(z))=R(q)$. Hence
\begin{equation}
y\left(\frac{1}{2\pi i}\log\left(\frac{\frac{4\pi}{i\sqrt[3]{4}}B_0\left(m^{*}(2z)^2;\frac{1}{6},\frac{2}{3}\right)}{f\left(\frac{4\pi}{i\sqrt[3]{4}}B_0\left(m^{*}(2z)^2;\frac{1}{6},\frac{2}{3}\right)\right)}\right)\right)=R(q).
\end{equation}
From (239) and (37) we get $G(R(q))=-P(Y(z))$. Combining this equation allong with (281) and (284) we get
$$
\frac{1}{P(Y(z))}=-\frac{1}{G(R(q))}=\sum^{\infty}_{n=1}nc_n\left(\frac{\frac{4\pi}{i\sqrt[3]{4}}B_0\left(m^{*}(2z)^2;\frac{1}{6},\frac{2}{3}\right)}{f\left(\frac{4\pi}{i\sqrt[3]{4}}B_0\left(m^{*}(2z)^2;\frac{1}{6},\frac{2}{3}\right)\right)}\right)^n.
$$
But it holds the next\\
\\
\textbf{Theorem 25.2} (Carty's theorem).\\
Whenever $Im(w)>0$ and $Im(z)>0$, we have
\begin{equation}
2\pi i\int^{w}_{z}\eta(t)^4dt=\left[\frac{1}{\sqrt[3]{4}}B_0\left(m^{*}(2t)^2;\frac{1}{6},\frac{2}{3}\right)\right]^{t=w}_{t=z},
\end{equation}
where $B_0(z;a,b):=\int^{z}_{0}t^{a-1}(1-t)^{b-1}dt$ is the incomplete Beta function.\\ 
\\

Hence if we define the function $m(z)$ such that
\begin{equation}
2\pi i \int^{m(z)}_{+i \infty}\eta(t)^4dt=z,
\end{equation}
then
\begin{equation}
\frac{1}{\sqrt[3]{4}}B_0\left(m^{*}(2m(z))^2;\frac{1}{6},\frac{2}{3}\right)=z
\end{equation}
and thus
$$
-\frac{1}{G(R(e(m(z))))}=\sum^{\infty}_{n=1}nc_n\left(\frac{\frac{4\pi}{i\sqrt[3]{4}}B_0\left(m^{*}(2m(z))^2;\frac{1}{6},\frac{2}{3}\right)}{f\left(\frac{4\pi}{i\sqrt[3]{4}}B_0\left(m^{*}(2m(z))^2;\frac{1}{6},\frac{2}{3}\right)\right)}\right)^n\Leftrightarrow
$$
$$
-\frac{1}{G(R(e(m(z))))}=\sum^{\infty}_{n=1}nc_n\left(\frac{\frac{4\pi}{i} z}{f\left(\frac{4\pi}{i}z\right)}\right)^n\Leftrightarrow
$$
$$
\frac{1}{\frac{4\pi z}{i}}-\frac{f'(\frac{4\pi z}{i})}{f(\frac{4\pi z}{i})}=-G(R(e(m(z)))).
$$
Hence we get the next\\
\\
\textbf{Theorem 25.3}\\
The quantities $f$ and $G$ and $P$ are related with the equation
\begin{equation}
-G(R(e(m(z))))=\frac{1}{z_1}-\frac{f'(z_1)}{f(z_1)}=P\left(\frac{1}{2\pi i}\log\left(\frac{z_1}{f\left(z_1\right)}\right)\right),
\end{equation}
where $z_1=\frac{4\pi z}{i}$.\\
\\
\textbf{Corollary 3.}\\
Knowing $G$ we can find $f$ from (291) and the oposite. The solution of the problem (233) is given by
\begin{equation}
y\left(\frac{1}{2\pi i}\log\left(\frac{\frac{4\pi}{i\sqrt[3]{4}}B_0\left(m^{*}(2z)^2;\frac{1}{6},\frac{2}{3}\right)}{f\left(\frac{4\pi}{i\sqrt[3]{4}}B_0\left(m^{*}(2z)^2;\frac{1}{6},\frac{2}{3}\right)\right)}\right)\right)=R(q),
\end{equation}
where $q=e(z)$, $Im(z)>0$. The function $Y(z)$ is given from
\begin{equation}
e(Y(z))=\frac{\frac{4\pi}{i\sqrt[3]{4}}B_0\left(m^{*}(2z)^2;\frac{1}{6},\frac{2}{3}\right)}{f\left(\frac{4\pi}{i\sqrt[3]{4}}B_0\left(m^{*}(2z)^2;\frac{1}{6},\frac{2}{3}\right)\right)}.
\end{equation}
\\
\textbf{Corollary 4.}\\
There exists constant $c\in\textbf{R}$, such that
\begin{equation}
2\pi i\int^{\frac{1}{2\pi i}\log\left(\frac{A}{f(A)}\right)}_{c}\frac{dt}{P(t)}=A.
\end{equation}
\\
\textbf{Proof.}\\
Integrate the second equality of (291).\\
\\
\textbf{Corollary 5.}\\
If $Im(z)>0$, then
\begin{equation}
-2\pi i \sqrt[3]{2}B_0\left(s\left(\frac{\log\left(\frac{A}{f(A)}\right)}{2\pi i}\right);\frac{1}{6},\frac{2}{3}\right)=A+c.
\end{equation}
\\
\textbf{Theorem 25.4}\\
If $f_0(A)=\sum^{\infty}_{n=1}\frac{f^{(n)}_0(0)}{n!}A^n$ and $a_n=\frac{f^{(n)}_0(0)}{\Gamma(n)}$, then
$$
Y^{(-1)}(z)=m\left(\frac{f_0(e(z))}{4\pi i}-\frac{c_1}{4\pi i}\right),
$$
$$
y(z)=R\left(e\left(m\left(\frac{f_0(e(z))}{4\pi i}-\frac{c_1}{4\pi i}\right)\right)\right).
$$
\\
\textbf{Proof.}\\
We have
$$
\frac{1}{P(z)}=\sum^{\infty}_{n=1}\frac{f^{(n)}_0(0)}{\Gamma(n)}e(nz)\Rightarrow \int^{Y(z)}_{Y(i\infty)}\frac{dt}{P(t)}=\frac{1}{2\pi i}\left[f_0\left(e\left(Y(z)\right)\right)-f_0\left(e(Y(i\infty))\right)\right]=
$$
$$
=\frac{-2}{\sqrt[3]{4}}B_0\left(m^*(2z)^2;\frac{1}{6},\frac{2}{3}\right).
$$
Hence
$$
\frac{f_0(e(Y(z)))}{4\pi i}-\frac{c_1}{4\pi i}=\frac{-1}{\sqrt[3]{4}}B_0\left(m^*(2z)^2;\frac{1}{6},\frac{2}{2}\right)\Rightarrow
$$
$$
m\left(\frac{f_0(e(Y(z)))}{4\pi i}-\frac{c_1}{4\pi i}\right)=z\Rightarrow Y^{(-1)}(z)=m\left(\frac{f_0(e(z))}{4\pi i}-\frac{c_1}{4\pi i}\right).
$$
\\

We are interested for products such
$$
\prod^{\infty}_{n=1}\left(1-q^n\right)^{-1/n\sum_{d|n}a_d\mu(n/d)}.
$$
Hence if we set $a_n=\frac{f_0^{(n)}(0)}{\Gamma(n)}$, then
$$
w(A)=f_0(A)
$$
and
$$
f(A)=\frac{A}{w^{(-1)}(A)}=\frac{A}{f_0^{(-1)}(A)}.
$$
Also then
$$
\prod^{\infty}_{n=1}\left(1-q_A^n\right)^{-1/n\sum_{d|n}a_d\mu(n/d)}=e^{f_0(A)}=e^{w(A)}\textrm{, }q_A=e(A)
$$
and
$$
\left[\left(\frac{D}{Dh}\right)^{n-1}(f(h))^{n}\right]_{h=0}=f_0^{(n)}(0).
$$
Hence
$$
y(A)=F_1\left(\frac{c}{2\pi i}-\frac{f_0(q_A)}{2\pi i}\right).
$$
Also if $g_1(A)$ is any function (analytic at the origin), with $g_1(0)=1$, we have
$$
g_1(f_0(A))=1+\sum^{\infty}_{n=1}\frac{A^n}{n!}\left[\left(\frac{D}{Dh}\right)^{n-1}(g_1'(h)f(h)^n)\right]_{h=0}.
$$
A simple consequence of the above is the following formula
$$
A\frac{w'(A)}{w(A)}=\sum^{\infty}_{n=0}\frac{A^n}{n!}\left[\left(\frac{D}{Dh}\right)^n(f(h))^n\right]_{h=0}.
$$
\\
\textbf{Theorem 25.5}\\
Assume the function $G(A)$. When $Im(A)>0$, $G(A)$ is given from the relation (i.e. satisfies the relation):
$$
4\pi i \eta(A)^4 G\left(R\left(e^{2\pi i A}\right)\right)=\psi(A).\eqno{(295.1)}
$$
Then $Y(A)=\int\psi(A)dA$ and if this integral is known (for example $\psi(A)$ is a polynomial), then $y(A)$ can be found explicity from $y\left(\int \psi(A)dA\right)=R\left(e^{2\pi i A}\right)$. Knowing $y(A)$ we can find $w\left(e^{2\pi i A}\right)$ from (see relation (29)): 
$$
w(q_A)=c-10\pi i \int^{y(A)}_{0}\frac{dt}{t\sqrt[6]{t^{-5}-11-t^5}}\textrm{, }q=e^{2\pi i A}.\eqno{(295.2)}
$$
Hence we find $f(A)$ from $f(A)=A/w^{(-1)}(A)$ and by this way evaluate all the above mentioned functions ($m_G(A)$, $Q(A)$, $P(A)$, $h(A)$, $s(A)$, $X(A),\ldots etc$) using only simple inversions and not using integral or series calculus.\\
\\
\textbf{Proof.}\\
From the relations $Y'(A)=-4\pi i \eta(A)^4P(Y(A))$ and $y(Y(A))=R(q)$, $G(y(A))=-P(A)$, we find 
$$
G(y(Y(A)))=-P(Y(A))\Leftrightarrow G(R(q))=-P(Y(A)).
$$
Hence $Y'(A)=4\pi i \eta(A)^4G(R(q))$. Hence given $G(A)$ such that $4\pi i \eta(A)^4G(R(q))=\psi(A)$ and $\int \psi(A)dA=\textrm{known}$, then $Y(A)$ is known. From this last fact we can evaluate all functions with simple inversions and without the use of calculus (integrals and infinite series).\\
\\
\textbf{Theorem 25.6}\\
There exists constant $c$ such that for all $A$ with $Im(A)>0$, we have
$$
Y(A)=2^{4/3}\int^{m^{*}(2A)}_{c}\frac{G(F_0(t))}{\left(t\sqrt{1-t^2}\right)^{2/3}}dt.\eqno{(295.3)}
$$
\\
\textbf{Proof.}\\
If we denote $\widetilde{m}^{*}(A)=\sqrt{1-m^{*}(A)^2}$, then we have
$$
4\pi i \eta(A)^4 G\left(R\left(e^{2\pi i A}\right)\right)=Y'(A)\Leftrightarrow 
$$
$$
2^{7/3}\frac{m^{*}{'}(2A)}{\left(\widetilde{m}^{*}(2A)m^{*}(2A)\right)^{2/3}}G\left(R\left(e^{2\pi i A}\right)\right)=Y'(A)\Leftrightarrow
$$
$$
2^{4/3}\frac{m^{*}{'}(A)}{\left(m^{*}(A)\sqrt{1-m^{*}(A)^2}\right)^{2/3}}G\left(F_0\left(m^{*}(A)\right)\right)=\frac{1}{2}Y'\left(\frac{A}{2}\right)\Leftrightarrow
$$
$$
Y\left(\frac{A}{2}\right)=2^{4/3}\int^{m^{*}(A)}_{c}\frac{G(F_0(t))}{\left(t\sqrt{1-t^2}\right)^{2/3}}dt.
$$
\\
\textbf{Theorem 25.7}\\
Assume that $G(A)$ is of the form 
$$
G\left(F_1(A)\right)=-\frac{1}{c-2\pi i A}+\psi_1'(c-2\pi i A)\textrm{, }Im(A)>0,\eqno{(295.4)}
$$
where $\psi_1(A)$ is given function, then we can find $w\left(e^{2\pi i A}\right)$ from Theorem 2 and hence solve the inversion problem 
$$
5\int^{y(A)}_{0}\frac{G(t)}{t\sqrt[6]{t^{-5}-11-t^5}}dt=A.
$$
Also for such $G(A)$ we can give evaluations to (119), using only inversion of functions and not calculus.\\
\\
\textbf{Remark.}\\
Actually if (295.4) holds then, we know $h(A)=\frac{1}{2\pi i}\left(\log(c-2\pi i A)-\psi_1(c-2\pi i A)\right)+c_1$, since we have $G(F_1(A))=h'(A)=-\frac{1}{c-2\pi i A}+\psi_1'(c-2\pi i A)$. Also then 
$$
f(A)=\exp\left(-2\pi i c_1+\psi_1(A)\right).\eqno{(295.5)}
$$
\\
\textbf{Proof.}\\
The proof is a immediate consequence of Theorem 2.\\
\\
\textbf{Example.}\\
Assume that $G$ is such that $G(F_1(A))=Ae^{A}$, then 
$$
5\int^{A}_{0}\frac{G(t)}{t\sqrt[6]{t^{-5}-11-t^5}}dt=
$$
$$
=5\int^{F_1^{(-1)}(A)}_{0}\frac{G(F_1(t))}{F_1(t)\sqrt[6]{F_1(t)^{-5}-11-F_1(t)^5}}F_1'(t)dt=
$$
$$
=\int^{F_1^{(-1)}(A)}_{0}G(F_1(t))dt=\int^{F_1^{(-1)}(A)}_{0}te^{t}dt=
$$
$$
=(F_1^{(-1)}(A)-1)e^{F_1^{(-1)}(A)}+1\Rightarrow
$$
$$
A=(F_1^{(-1)}(y(A))-1)e^{F_1^{(-1)}(y(A))}+1\Rightarrow
$$
$$
y_i(F_1(A))=(A-1)e^A+1\Rightarrow y(A)=F_1\left[\left((t-1)e^t+1\right)^{(-1)}(A)\right].
$$
Hence if we denote $\psi_1(t)=(t-1)e^t+1$, then $\psi_1'(t)=te^t$ and
$$
y(A)=F_1(\psi_1^{(-1)}(A)).
$$

\section{Constructing modular forms}

If $f(z)$ is analytic and have no roots in $\textbf{C}$ and $g_{n}(z)=(z-z_1)(z-z_2)\ldots (z-z_n)$, we define the function
$$
F(z)=\frac{f(z)}{g_n(z)}\textrm{, }z\in\textbf{C}.
$$ 
Then if $z_k(w)=z^{(0)}_{k}(w)$, are the roots of $\frac{1}{F(z)}=w$, with $z^{(0)}_{k}(0)=z_k$, $k=1,2,\ldots,n$ and $z^{(j)}_k(0)=\left(\frac{d^j}{dw^j}z_k(w)\right)_{w=0}$, for all $j=1,2,\ldots$, then
$$
\oint_{C}h(F(z))dz=2\pi i\sum^{n}_{k=1}\sum^{\infty}_{j=1}\frac{h^{(j)}(0)}{\Gamma(j)}\frac{z^{(j)}_k(0)}{j!}.
$$
The curve $C$ is any simple closed curve that enclose $z_1,z_2,\ldots,z_n$.\\
In case we have $F(z)=\frac{f(z)}{g(z)}$, with $f(z)$ analytic and $f(z)\neq 0$ in a neighborhood of $0$ and $g(z)$ analytic in a neighborhood of $0$ with simple zero at $z=0$, then
\begin{equation}
\oint_{C}h(F(z))dz=2\pi i \sum^{\infty}_{k=1}\frac{h^{(k)}(0)}{\Gamma(k)}\frac{1}{k!} \left(\left(\frac{D}{Dt}\right)^{k-1}\left(tF(t)\right)^{k}\right)_{t=0}.
\end{equation}
\\
\textbf{Theorem 26.}\\
If $f(z)$ is entire function with no poles and zeroes in $\textbf{C}$ and $g(z)$ entire with $g(0)\neq 0$ and $z_1,z_2,\ldots$ all the zeros (simple) such that $\sum^{\infty}_{n=1}\frac{1}{|z_n|}<\infty$, then exists constants $a,b\in\textbf{C}$ such that 
\begin{equation}
g(z)=e^{az+b}\prod^{\infty}_{n=1}\left(1-\frac{z}{z_n}\right).
\end{equation}
Also if $h(z)$ also analytic in $\textbf{C}$, then for $F(z)=\frac{f(z)}{g(z)}$ holds
\begin{equation}
\oint_{C}h(AF(z))dz=2\pi i \sum^{\infty}_{j,k=1}\frac{h^{(k)}(0)}{\Gamma(k)}\frac{A^k}{k!} \left(\left(\frac{D}{Dt}\right)^{k-1}\left(\left(t-z_j\right)F(t)\right)^{k}\right)_{t=z_j},
\end{equation}
where $C$ is a simple closed path that enclose all $z_1,z_2,\ldots$\\
\\
\textbf{Theorem 26.1}\\
Assume the function $f(x)$ with $f(x)=g(x)/x$, $g(x)$ being analytic arround 0 with radius at least 1 and $g(0)\neq 0$. If
$$
\int^{2\pi}_{0}\log\left(1-Af(e^{it})\right)e^{it}dt=-2\pi \rho(A),\eqno{(298.1)}
$$
then for every $A$ with $|A|\leq |A_0|$, $\rho(A)$ is solution of the equation $1/f(x)=A$, where $A_0$ is a certain complex number.\\
\textbf{Remarks.}\\
\textbf{i)} Note that (298.1) is a general integral that can be found only by inversion of function. That is, by solving an equation and does not requires the knowledge of derivatives or other integrals.\\   
\textbf{ii)} From relation (298.1) and assuming that $f$ satisfies the requirements of the theorem, we can easily get
$$
\left|\rho(A)\right|\leq \left|\log\left(1-M_f\cdot |A|\right)\right|\textrm{, }\forall A \in D_f,\eqno{(298.2)}
$$
where 
$$
D_{f}=\left\{A\in\textbf{C}:|A|\leq |A_0|< 1/M_f\right\}\eqno{(298.3)}
$$ 
and 
$$
M_f=\textrm{sup}\left\{\left|f\left(e^{it}\right)\right|:t\in[0,2\pi)\textrm{, }f(A) \textrm{ analytic in }|A|\leq 1\right\}.\eqno{(298.4)}
$$
\\ 
\textbf{Proof.}\\
Assume the equation $1/f(x)=A:(eq)$, $f(x)=g(x)/x$, with $g(x)$ analytic arround 0 and $g(0)\neq 0$. Then a solution of $(eq)$ is  
$$
x=\sum^{\infty}_{n=1}\frac{A^n}{n!}\left[\left(\frac{D}{Dh}\right)^{n-1}(g(h))^n\right]_{h=0}.
$$ 
Also from Theorem 26 we can extract easily the identity
$$
\frac{1}{2\pi i}\oint_{C}\left(f(z)\right)^ndz=\frac{1}{\Gamma(n)}\left[\left(\frac{D}{Dh}\right)^{n-1}(hf(h))^n\right]_{h=0},
$$
where $C$ is a simple closed path enclouses 0 i.e. we can select $|z|=r\leq 1$. Hence
$$
\frac{1}{2\pi}\int^{2\pi}_{0}\left(f(e^{it})\right)^ne^{it}dt=\frac{1}{\Gamma(n)}\left[\left(\frac{D}{Dh}\right)^{n-1}(hf(h))^n\right]_{h=0}.
$$
Using the above identities we get easily the result.\\ 
\\

Assume now the function 
$$
f(z)=\exp\left(c_0+c_1z+c_2z^2+\ldots\right).
$$
Set
$$
g(z)=\left(\frac{i}{\sqrt[a]{d}}\sqrt[a]{\frac{f(1-z)}{f(z)}}\right)^{-1}.
$$
If $B=w(A)$ is root of the equation
\begin{equation}
g(B)=A,
\end{equation}
then the Lagrange inversion theorem states (in a different vcersion): If $g(A)$ is analytic arround $0$ and $g'(0)\neq 0$, equation (299), have solution 
$$
B=w(A)=\sum^{\infty}_{n=1}c_n (A-g(0))^n,
$$
where
\begin{equation}
c_n=\lim_{h\rightarrow 0}\frac{1}{n!}\left[\frac{D^{n-1}}{Dh^{n-1}}\left(\frac{h}{g(h)-g(0)}\right)^n\right].
\end{equation}
Also for the function $w(A)$ holds
$$
w\left(-\frac{1}{D_1A}\right)+w(A)=1,
$$
where $D_1=d^{-2/a}$. But then also
\begin{equation}
c_n=\frac{1}{2\pi i n}\oint_{C}\left(\frac{1}{g(z)-g(0)}\right)^ndz.
\end{equation}
Hence
$$
w(A)=\sum^{\infty}_{n=1}c_n(A-g(0))^n=\frac{1}{2\pi i}\sum^{\infty}_{n=1}\frac{(A-g(0))^n}{n}\oint_{C}\frac{1}{(g(z)-g(0))^n}dz.
$$
Hence
$$
w(A)=\frac{1}{2\pi i}\oint_{C}\log\left(\frac{A-g(z)}{g(0)-g(z)}\right)dz=\frac{1}{2\pi i}\oint_{C}\log\left(1-\frac{g(0)-A}{g(0)-g(z)}\right)dz.
$$
The function $w(A)$ is $T-$periodic iff
$$
\oint_{C}\log\left(1+\frac{T}{A-g(z)}\right)dz=0
$$
and we have the next\\
\\
\textbf{Theorem 27.}\\
Assume the function
\begin{equation}
f(z)=\exp\left(c_0+c_1z+c_2z^2+\ldots\right).
\end{equation}
Set
\begin{equation}
g(z)=\left(\frac{i}{\sqrt[a]{d}}\sqrt[a]{\frac{f(1-z)}{f(z)}}\right)^{-1}.
\end{equation}
If $B=w(A)$ is the solution of the equation
\begin{equation}
g(B)=A,
\end{equation}
then
\begin{equation}
w(A)=\sum^{\infty}_{n=1}c_n(A-g(0))^n,
\end{equation}
where
\begin{equation}
c_n=\lim_{h\rightarrow 0}\frac{1}{n!}\left[\frac{D^{n-1}}{Dh^{n-1}}\left(\frac{h}{g(h)-g(0)}\right)^n\right]=\frac{1}{2\pi i n}\oint_{C}\left(\frac{1}{g(z)-g(0)}\right)^ndz.
\end{equation}
Also for the function $w(A)$ holds
\begin{equation}
w\left(-\frac{1}{D_1A}\right)+w(A)=1,
\end{equation}
where $D_1=d^{-2/a}$ and
\begin{equation}
w(A)=\frac{1}{2\pi i}\oint_{C}\log\left(1-\frac{g(0)-A}{g(0)-g(z)}\right)dz.
\end{equation}
\\
\textbf{Theorem 28.}\\
For the function $H(A)=f(w(A))$, with $f,w$ as in the above theorem we have
\begin{equation}
H\left(-\frac{1}{D_1A}\right)=i^{-a} d A^{-a}H(A)\textrm{, }Im(A)>0
\end{equation}
and $D_1=d^{-2/a}$. Moreover $w(A)$ is $T-$periodic iff
\begin{equation}
\oint_{C}\log\left(1+\frac{T}{A-g(z)}\right)dz=0,
\end{equation}
for all $A$, with $Im(A)>0$, $0<Re(A)<\sigma$, where $\sigma$ is a certain positive real.\\
\textbf{Remark.} A condition under $w(A)$ is periodic, is 
\begin{equation}
Im(g(z))\leq 0.
\end{equation}

\section{Appendix: A Gauss hypergeometric function with some curious functional equations}

Let
\begin{equation}
F(x):={}_2F_1(a,b;2b;x),
\end{equation}
then
\begin{equation}
F(x)=(1-x)^{-a}F\left(\frac{x}{x-1}\right)
\end{equation}
Consider now the tranformation
$$
w=\frac{A x+B}{C x+D}
$$
with
$$
\frac{F\left(\frac{A y+B}{Cy+D}\right)}{F(y)}=n\frac{F\left(\frac{Ax+B}{Cx+D}\right)}{F(x)},
$$
where $y=\frac{x}{x-1}$, then a choice of $A,B,C,D$ is $D=-A$, $B=-2A$ and
$$
n=\left(\frac{(x-1)[(C+A)x-3A]}{A-Cx}\right)^{-a}.
$$
Hence we get the following\\
\\
\textbf{Theorem 1.}\\
Set $\lambda=\frac{C}{2A}$, then equation
$$
\frac{F\left(\frac{1-x/2}{1/2-\lambda x}\right)}{F(x)}=\sqrt{r},
$$
have root $\beta_r$ such that
$$
\beta_{n^2r}=\frac{\beta_r}{\beta_r-1},
$$
where
$$
n=\left(-\frac{(1-\beta_r)(1-2\lambda \beta_r)}{1+(2\lambda-1)\beta_r}\right)^{-a}.
$$
The multiplier is
$$
m=\frac{F(\beta_{n^2r)}}{F(\beta_r)}=(1-\beta_r)^a.
$$

From another point of view if we set $\beta^{(-1)}(x)$ to be the inverse function of $\beta_x=\beta(x)$, we have 
$$
\beta\left(x\left(-\frac{(1-\beta(x))(1-2\lambda\beta(x))}{1+(2\lambda-1)\beta(x)}\right)^{-2a}\right)=\frac{\beta(x)}{\beta(x)-1}.
$$
Hence
$$
\beta\left(\beta^{(-1)}(x)\left(-\frac{(1-x)(1-2\lambda x)}{1+(2\lambda-1)x}\right)^{-2a}\right)=\frac{x}{x-1}
$$
and hence
$$
\beta^{(-1)}(x)\left(-\frac{(1-x)(1-2\lambda x)}{1+(2\lambda-1) x}\right)^{-2a}=\beta^{(-1)}\left(\frac{x}{x-1}\right).\eqno{(f)}
$$
Also
$$
\beta^{(-1)}(x)=\left(\frac{F\left(\frac{1-x/2}{1/2-\lambda x}\right)}{F(x)}\right)^2.
$$
Hence for the evaluation of $\beta^{(-1)}(x)$, i.e. when we know the value of ceratin $x=\beta_r$ and we want to evaluate the values $r$, then we can use formula $(f)$, or we can use modular relations.\\
\\

Set now $G(z)=F(1-z^2)$. Then for every $z\in\textbf{C}^{*}$ we have
\begin{equation}
G\left(-\frac{1}{z}\right)=z^{2a}G(z).
\end{equation}
Now assume the functions
\begin{equation}
\theta_2(z)=\sum^{\infty}_{n=-\infty}q^{(n+1/2)^2}\textrm{, }q=e(z)
\end{equation}
and
\begin{equation}
\theta_3(z)=\sum^{\infty}_{n=-\infty}q^{n^2}\textrm{, }q=e(z),
\end{equation}
where $e(z):=e^{2\pi i z}$ and $Im(z)>0$. Then the function
\begin{equation}
m(z):=\left(\frac{\theta_2(z)}{\theta_3(z)}\right)^4,
\end{equation}
satisfies the relations
\begin{equation}
m(z+2)=m(z),
\end{equation}
\begin{equation}
m(z+1)=\frac{m(z)}{m(z)-1}
\end{equation}
and
\begin{equation}
m\left(-\frac{1}{z}\right)=1-m(z).
\end{equation}
Hence if we define the function
\begin{equation}
H(z):=F(m(z))={}_2F_1(a,b;2b;m(z)),
\end{equation}
then using the above properties, we have
\begin{equation}
H(z+1)=(1-m(z))^aH(z).
\end{equation}
Clearly 
\begin{equation}
H(z+2)=H(z)
\end{equation}
and
\begin{equation}
H\left(\frac{z}{z+1}\right)=F\left(\frac{1}{m(z)}\right).
\end{equation}
Set $\frac{z}{z+1}=w$. Then $z=\frac{-1}{w+1}$ and we have
$$
H\left(-\frac{1}{w}\right)=F\left(\frac{1}{m\left(-\frac{1}{w+1}\right)}\right)=F\left(\frac{1}{1-m\left(w+1\right)}\right)=F\left(\frac{1}{1-\frac{m(w)}{m(w)-1}}\right)=
$$
$$
=F\left(1-m(w)\right)=F\left(\frac{m(w)-1}{m(w)}\right)m(w)^{-a}=F\left(\frac{1}{m(w+1)}\right)m(w)^{-a}=
$$
$$
=H\left(\frac{w+1}{w+2}\right)m(w)^{-a}.
$$
Hence we get\\
\\
\textbf{Theorem 2.}
\begin{equation}
H\left(-\frac{1}{w}\right)=H\left(\frac{w+1}{w+2}\right)m(w)^{-a}\textrm{, }Im(w)>0.
\end{equation}
An obsevation is that if $a_{i,j}\in\textbf{Z}$, for $i=1,2,3,4$, $j=1,2$ with $a_{1,j}a_{4,j}-a_{2,j}a_{3,j}=1$, $j=1,2$. Then if $a_{i,1}\equiv a_{i,2} (\textrm{mod} 2)$, $i=1,2,3,4$, we have
\begin{equation}
H\left(\frac{a_{1,1}z+a_{2,1}}{a_{3,1}z+a_{4,1}}\right)=H\left(\frac{a_{1,2}z+a_{2,2}}{a_{3,2}z+a_{4,2}}\right)\textrm{, }Im(z)>0.
\end{equation} 
\\

Hence we can write relations like
\begin{equation}
H\left(-\frac{1}{z}\right)=H\left(-\frac{1}{z+2n}\right)\textrm{, }n\in\textbf{Z}
\end{equation}
and
\begin{equation}
H\left(-\frac{1}{z}\right)=H\left(2n+1-\frac{1}{z+2l}\right)m(z)^{-a}\textrm{, }n,l\in\textbf{Z}
\end{equation}
and
\begin{equation}
H\left(\frac{z}{2n z+1}\right)=H(z)\textrm{, }n\in\textbf{Z}.
\end{equation}

\[
\]

\centerline{\bf References}\vskip .2in

[1]: N.D. Bagis. ''Generalized Elliptic Integrals and Applications''. arXiv:1304.2315v2 [math.GM]. 22 Jun 2013\\ 

[2]: N.D. Bagis. ''On Generalized Integrals and Ramanujan-Jacobi Special Functions''. arXiv:1309.7247v5 [math.GM] 15 Nov 2021\\

[3]: N.D. Bagis. ''On the Complete Evaluation of Theta Functions''. arXiv:1503.01141v4 [math.GM] 10 Mar 2021\\

[4]: J.V. Armitage, W.F. Eberlein. ''Elliptic Functions''. Cambridge University Press. (2006)\\

[5]: N.D. Bagis. ''On the Construction of Relativistic Quantum Wave Equation and General Solution of the Second Order  Differential Equation''.\\

[6]: N.D. Bagis, M.L. Glasser. ''Conjectures on the evaluation of certain functions with algebraic properties''. Journal of Number Theory 155 (2015) 63-84 (Elsevier).\\

\[
\]

\end{document}